\documentclass[12pt]{amsart}

\usepackage{amsmath,amscd,amsfonts,amsthm,amssymb,eucal}
\usepackage{url}
\usepackage{fullpage}

\newcommand{\mods}[1]{\,(\mathrm{mod}\,{#1})}

\renewcommand{\leq}{\leqslant}
\renewcommand{\geq}{\geqslant}

\newcommand{\fleche}[1]{\stackrel{#1}{\longrightarrow}}
\newcommand{\eps}{\varepsilon}
\DeclareMathOperator{\SL}{SL}

\DeclareMathOperator{\GL}{GL}

\DeclareMathOperator{\End}{End}
\DeclareMathOperator{\Lie}{Lie}
\DeclareMathOperator{\Spec}{Spec}

\DeclareMathOperator{\Gal}{Gal}

\DeclareMathOperator{\Jac}{Jac}

\DeclareMathOperator{\Aut}{Aut}
\DeclareMathOperator{\Sp}{Sp}
\DeclareMathOperator{\GSp}{GSp}

\newcommand{\field}[1]{\mathbb{#1}}

\newcommand{\Q}{\field{Q}}

\newcommand{\Z}{\field{Z}}
\newcommand{\F}{\field{F}}
\newcommand{\Fp}{\field{F}_p}

\newcommand{\R}{\field{R}}
\newcommand{\C}{\field{C}}
\renewcommand{\P}{\field{P}}
\newcommand{\A}{\field{A}}

\newcommand{\EE} {\mathcal{E}}

\newcommand{\ra}{\to}

\newcommand{\OO}{\mathbb{Z}}

\newcommand{\inj}{\hookrightarrow}

\newcommand{\G}{\field{G}}

\newcommand{\beq}{\begin{displaymath}}
\newcommand{\eeq}{\end{displaymath}}
\newcommand{\beqn}{\begin{equation}}
\newcommand{\eeqn}{\end{equation}}

\newcommand{\Fbar}{\bar{\F}}

\newcommand{\lto}{\longrightarrow}

\newcommand{\barU}{\mathbf{U}}
\newcommand{\barG}{\mathbf{G}}
\newcommand{\barN}{\mathbf{N}}
\newcommand{\barT}{\mathbf{T}}
\newcommand{\barGL}{\mathrm{\bf GL}}
\newcommand{\barI}{\mathbf{I}}
\newcommand{\barSL}{\mathrm{\bf SL}}
\newcommand{\g}{\mathfrak{g}}
\newcommand{\M}{\mathrm{M}}
\newcommand{\gp}[1]{\langle #1 \rangle}
\newcommand{\w}{\mathfrak{w}}

\newcommand{\et}{\mathrm{et}}

\newcommand{\opcit}{\textit{op.~cit.}}

\newcounter{ccounter}
\renewcommand{\c}{\addtocounter{ccounter}{1}c_{\arabic{ccounter}}}

\theoremstyle{plain}
\newtheorem{thm}{Theorem}
\newtheorem{prop}[thm]{Proposition}
\newtheorem{cor}[thm]{Corollary}
\newtheorem{lem}[thm]{Lemma}

\theoremstyle{definition}
\newtheorem{defn}[thm]{Definition}

\newtheorem{exmp}[thm]{Example}
\newtheorem{que}[thm]{Question}

\theoremstyle{remark}
\newtheorem{rem}{Remark}

\begin{document}
\title[Expanders and Galois representations]{Expander graphs,
  gonality and variation of Galois representations}

\author[J. Ellenberg]{Jordan S. Ellenberg}
\address{Department of Mathematics \\ 
University of Wisconsin \\ 480 Lincoln Drive \\
Madison, WI 53705 USA} 
\email{ellenber@math.wisc.edu}

\author[C. Hall]{Chris Hall}
\address{Department of Mathematics
University of Wyoming, Ross Hall,
Laramie, WY 82071, USA} \email{chall14@uwyo.edu}

\author[E. Kowalski]{Emmanuel Kowalski}
\address{ETH Z\"urich -- D-MATH\\
  R\"amistrasse 101\\
  8092 Z\"urich\\
  Switzerland} \email{kowalski@math.ethz.ch}



\begin{abstract}
  We show that families of coverings of an algebraic curve where the
  associated Cayley-Schreier graphs form an expander family exhibit
  strong forms of geometric growth. We then give many arithmetic
  applications of this general result, obtained by combining it with
  finiteness statements for rational points of curves with large
  gonality.  In particular, we derive a number of results concerning
  the variation of Galois representations in one-parameter families of
  abelian varieties. 
\end{abstract}

\maketitle

%\tableofcontents

\section{Introduction}

When $A \ra B$ is a family of abelian varieties over a base $B$, there
is a general philosophy that "most" fibers $A_b$, where $b$ ranges
over closed points of $B$, should have properties similar to that of
the generic fiber $A_\eta$.  In the present paper we develop a very
general method to prove statements of this kind, where the properties
of abelian varieties we study pertain to the images of their
$\ell$-adic Galois representation, or to the presence of ``extra''
algebraic cycles, and where the base $B$ is a curve over a number
field.  Notably, a key role is played by recent results about
expansion in Cayley graphs of linear groups over finite fields.

Our motivation for this work comes from our previous
paper~\cite{eehk}, joint with C. Elsholtz, in which we studied the
geometrically non-simple specializations in a family of abelian
varieties whose generic fiber is geometrically simple.  As pointed out
by J.~Achter (and indicated in a note of~\cite{eehk}), D.~Masser
previously studied questions of a similar nature in \cite{masser}
using methods arising from transcendence theory.
\par
A particular case of our results and of those of Masser is the
following: let $g\geq 2$, let $f\in \Z[X]$ be a squarefree
polynomial of degree $2g$, and consider the family of hyperelliptic
curves
\begin{equation}\label{eq:concrete}
  \mathcal{C}_t\,:\,y^2=f(x)(x-t),\quad\quad 
  t\in U=\A^1-\{\text{zeros of $f$}\}.
\end{equation}
\par
The jacobian $J_t=\Jac(\mathcal{C}_t)$ is an abelian variety of
dimension $g$.  In fact, the jacobian of the generic fiber,
$\Jac(\mathcal{C}_{\eta})$, is an absolutely simple abelian variety
with geometric endomorphism ring equal to $\Z$.  Masser's methods, as
well as ours, can be interpreted as stating that, for ``most'' $t$,
the specialization $J_t$ is absolutely simple with geometric
endomorphism ring $\Z$.  In \cite{eehk}, the parameter $t$ varies over
$\Q$; by contrast, Masser's methods provide similar results where $t$
is allowed to range over the union of all extensions of degree $d$ of
$\Q$, for any fixed $d \geq 1$.
%However, whereas we could
%handle points $t\in k$ for a fixed number field $k$, he was able to
%handle simultaneously all numbers fields $k$ of degree $\leq d$, for
%$d\geq 1$, and give upper bounds for those of bounded height for which
%the endomorphism ring of $J_t$ is not $\Z$.
%  We then want to extend this property to
% ``most'' fibers, and deduce as corollaries that there is a similar
% persistence for the endomorphism ring and for geometric simplicity.
\par
The prototype of a statement that involves all number fields of degree
$d$ is the deep fact that there exist only finitely many
$j$-invariants of CM elliptic curves with bounded degree (or
equivalently, in view of the theory of complex multiplication, the
class number of imaginary quadratic orders tends to infinity with the
absolute value of the discriminant). Our method allows us to prove
many new finiteness statements of this nature.

% In the present paper, we present
% a new method to prove strong forms of persistence of generic
% properties for variations of Galois representations which, like
% Masser's, applies to fibers over number fields of bounded degree, not
% only over a fixed number field.

\subsection{Statements of results}

The new method is quite general. It is based, rather surprinsingly, on
expansion properties of some families of graphs. We encapsulate a
sufficiently strong form of expansion in the following definition:
% Although the best
% known examples are expander graphs, our results are applicable for
% families satisfying a weaker assumption, which we call
% \emph{esperantism}.

\begin{defn}[Esperantist graphs]\label{def:esperantist}
  Let $(\Gamma_i)$ be a family of connected $r$-regular
  graphs,\footnote{\ We allow our graphs to have self-loops and
    multiple edges between vertices; see, e.g.,~\cite[\S
    4.2]{lubotzky} for references.} with adjacency matrices
  $A(\Gamma_i)$. Let
$$
\Delta_i=r\mathrm{Id}-A(\Gamma_i)
$$
denote the combinatorial Laplace operator of $\Gamma_i$, and let 
$$
0=\lambda_0(\Gamma_i)<\lambda_1(\Gamma_i)\leq \lambda_2(\Gamma_i)\leq
\ldots
$$
denote the spectrum of $\Delta_i$. We say that $(\Gamma_i)$ is an
\emph{esperantist family} if it satisfies\footnote{\ This limit means
  that for any $N\geq 1$, there are only finitely many $i\in I$ such
  that $|\Gamma_i|\leq N$; similar limits below have analogue
  meaning.}
\begin{equation}\label{eq:limit-size}
\lim_{i\ra +\infty}{|\Gamma_i|}=+\infty,
\end{equation}
and if there exist constants $c>0$ and $A\geq 0$, such that
\begin{equation}\label{eq:esperanto}
\lambda_1(\Gamma_i)\geq \frac{c}{(\log 2|\Gamma_i|)^A}
\end{equation}
for all $i$.
\end{defn}

\begin{exmp}
  In the special case where we can take $A=0$, we obtain the
  well-known notion of an \emph{expander family} of graphs, namely
\begin{equation}\label{eq:def-expand}
\lambda_1(\Gamma_i)\geq c>0
\end{equation}
for all $i$ and some constant $c$ independent of $i$. We refer readers
to the survey paper of Hoory, Linial and Wigderson~\cite{hlw} for
extensive background information on expansion in graphs, in particular
on expanders. In Section~\ref{sec-beyond}, we will comment on the
distinction we make between expander graphs and esperantist graphs.
% \par
% In fact, one could weaken further the expansion requirement, at least
% for the applications contained in this paper, see
% Section~\ref{sec-beyond} at the end of the paper.
\end{exmp}

The following theorem is our main diophantine statement:

\begin{thm}\label{th:main}
  Let $k$ be a number field.  Let $U/k$ be a smooth geometrically
  connected algebraic curve over $k$ and $(U_i)_{i\in I}$ an infinite
  family of \'etale covers of $U$ defined over $k$.  Fix an embedding
  of $k$ in $\C$ to define the complex Riemann surfaces $U_\C$ and
  $U_{i,\C}$. Let $S$ be a fixed finite symmetric\footnote{\ This
    means that $s^{-1}\in S$ for all $s\in S$.} generating set of the
  topological fundamental group $\pi_1(U_{\C},x_0)$ for some fixed
  $x_0\in U$ and assume that the family of Cayley-Schreier graphs
  $C(N_i,S)$ associated to the finite quotient sets
$$
N_i=\pi_1(U_{\C},x_0)/\pi_1(U_{i,\C},x_i),\quad\quad x_i\in
U_{i}\text{ some point over } x_0,
$$
is an \emph{esperantist} family.
\par
Then, for any fixed $d\geq 1$, the set
$$
\bigcup_{[k_1:k]=d}{U_i(k_1)}
$$
is finite for all but finitely many $i$.
\end{thm}

% \begin{rem} 
%   We denote by $U_\C$ the complex manifold formed by the complex
%   points of $U$, after the choice of some embedding $k \inj \C$.
% \end{rem}

\begin{rem}
  We recall the definition of a Cayley-Schreier graph: given a group
  $G$ with symmetric generating set $S\subset G$ and a subgroup $H$,
  the graph $C(G/H,S)$ is defined as the $|S|$-regular (undirected)
  graph with vertex set $G/H$ and with (possibly multiple) edges from
  $xH$ to $sxH$ for all $s\in S$. If $H=1$, this is the Cayley graph
  of $G$ with respect to $S$.
 % (We assume that our index set $I$ is infinite, since our claims
 %  are vacuous when it is finite, but it makes sense to state that a
 %  finite collection of graphs is an expander).
\par
Theorem~\ref{th:main} applies, in particular, when the family is an
expander family, which is the case in many of the applications in this
paper.
\par
% As we will explain later, all the graphs occurring in this paper are
% very likely to be expander families, and this has been established
% in many cases. However, the esperantist property may be quite a bit
% easier to check, and the most recent developments proving expansion
% often begin with a proof of the weaker esperantist
% property. (Exceptions are arguments based on Property $(T)$ or on
% automorphic methods.)
% \par
As far as we know, Theorem~\ref{th:main} (and its variants and
applications) are the first explicit use of general theorems about
spectral gaps in graphs to obtain finiteness statements in arithmetic
geometry.  However the idea descends from a result of Zograf~\cite{z}
(and, independently, Abramovich~\cite{a}), who proved lower bounds for
gonality of modular curves via spectral gaps for the Laplacian on the
underlying Riemann surfaces (given by Selberg's $3/16$ bound).
\par
We believe that there should exist other arithmetic consequences of
the following philosophy: in any category where the notion of finite
(Galois) covering makes sense, any family of finite coverings with a
similar expansion property should be ``extremely complicated''. In
Section~\ref{sec:conclusion}, we present some concrete questions along
these lines. 
\end{rem}

In the remainder of the introduction, we present some applications of
this theorem. The reader may note that these do not mention explicitly
a family of coverings of a fixed curve: those appear only as auxiliary
tools in the proofs. The finiteness statements are therefore of a
different nature than the finiteness of the set of rational points on
some fixed algebraic variety. We focus on abelian varieties, but this
is by no means the only area where applications are possible.
\par
% In Section~\ref{sec:application}, we start with some fairly direct
% applications, working with abelian varieties of specific type.
First, we will give a strong uniform version of ``large Galois image''
for certain one-parameter families of abelian varieties:

\begin{thm}\label{th:sp}
  Let $k$ be a number field and $U/k$ a smooth geometrically connected
  algebraic curve over $k$.  Let $\mathcal{A}\ra U$ be a principally
  polarized abelian scheme of dimension $g\geq 1$, defined over $k$,
  and let
$$
\rho\,:\, \pi_1(U_{\C},x_0)\ra \Sp_{2g}(\Z)
$$
be the associated monodromy representation.  For $k_1/k$ a finite
extension and $t\in U(k_1)$, let $\bar{\rho}_{t,\ell}$ be the Galois
representation
$$
\bar{\rho}_{t,\ell}\,:\, \Gal(\bar{k}/k_1)\ra \GSp_{2g}(\F_{\ell})
$$
associated to the action on the $\ell$-torsion points of
$\mathcal{A}_t$.  
\par
If the image of $\rho$ is Zariski-dense in $\Sp_{2g}$, then for any
$d\geq 1$ and all but finitely many $\ell$, depending on $d$, the set
$$
\bigcup_{[k_1:k]=d}{\{t\in U(k_1)\,\mid\, \text{the image of
    $\bar{\rho}_{t,\ell}$ does not contain $\Sp_{2g}(\F_{\ell})$} \}}
$$
is finite.
\end{thm}

We proved the case $d=1$ in our earlier paper~\cite[Prop. 8]{eehk}.

\begin{cor}\label{cor:concrete}
Let $k$ be a number field, and let $f\in k[X]$ be a squarefree polynomial
of degree $2g$ with $g\geq 1$. Let $U_f$ be the complement of the
zeros of $f$ in $\A^1$, and let $\mathcal{C}/U$ be the family of
hyperelliptic curves given by
$$
\mathcal{C}\,:\, y^2=f(x)(x-t),
$$ 
with Jacobians $J_t=\Jac(\mathcal{C}_t)$. Then for any $d\geq 1$, the
set
$$
\bigcup_{[k_1:k]=d}{
\{t\in U(k_1)\,\mid\, \End_{\C}(J_t)\not=\Z\}
}
$$
is finite.
\end{cor}

\begin{rem}
  This set was also considered by Masser, who gave an explicit upper
  bound for the cardinality of its subsets of bounded height
  (see~\cite[Theorem, p. 459]{masser}), namely we have
$$
\Bigl|
\bigcup_{[k_1:k]=d}{
\{t\in U(k_1)\,\mid\, \End_{\C}(J_t)\not=\Z\text{ and } 
h(t)\leq h\}
}
\Bigr|
\ll \max(g,h)^{\beta}
$$
for $h\geq 1$, where $h(x)$ denotes the absolute logarithmic Weil
height on $\bar{k}$ and $\beta>0$ is some (explicit) constant
depending only on $g$. Then, based on concrete examples and other
results of Andr\'e, Masser raised the following question
(see~\cite[middle of p. 460]{masser}): is it true, or not, that there
are only finitely many $t$ of degree at most $d$ over $k$ such that
the geometric endomorphism ring of $J_t$ is larger than $\Z$?  This
corollary gives an affirmative answer.
\par
That being said, contrary to Masser's method, ours does not give
explicit or effective bounds on the cardinality of the sets we
consider, and hence the two are complementary.  (We discussed a
similar dialectic in ~\cite{eehk}.)  In particular, Masser's methods
can be used to get some control of exceptional fibers in families of
abelian varieties over higher-dimensional bases, whereas
Theorem~\ref{th:main} does not provide anything interesting in such a
situation.
\end{rem}

We will next prove that Theorem~\ref{th:main} implies that two
families of elliptic curves which are not generically isogenous have
few fibers with isomorphic mod-$\ell$ Galois representations:

\begin{thm}\label{th:sl2sl2}
  Let $k$ be a number field and let $\EE_1$ and $\EE_2$ be
  elliptic curves over the function field $k(T)$.
  Assume furthermore that $\EE_1$ and $\EE_2$ are not geometrically
  isogenous.  Then, for $d \geq 1$, the set
$$
\bigcup_{[k_1:k]=d}{\{t\in k_1\,\mid\, \text{$\EE_{1,t}[\ell]$ and
    $\EE_{2,t}[\ell]$ are isomorphic as $G_{k_1}$-modules} \}}
$$
is finite for all but finitely many $\ell$.
\end{thm}

Our method has also some applications to \emph{arbitrary}
one-parameter families of abelian varieties. Using a general
``semisimple approximation'' of the Galois groups of $\ell$-torsion
fields (which builds on work of Serre~\cite{serre:vig}), we will prove
the following result on existence, which also exhibits some level of
quantitative information on the field of definition (the interest of
the latter was suggested by J. Pila):

\begin{thm}\label{th:torsion}
  Let $k$ be a number field and let $U/k$ a smooth geometrically
  connected algebraic curve over $k$.  Let $\mathcal{A}\ra U$ be an
  abelian scheme of dimension $g\geq 1$, defined over $k$.  Then for every
  $d \geq 1$ there exists an $\ell(d)$ such that, for all primes $\ell
  > \ell(d)$, the set
\begin{equation}\label{eq-exist-torsion}
\bigcup_{[k_1:k]=d}{\{t\in U(k_1)\,\mid\,
  \text{$\mathcal{A}_t[\ell](k_1)$ is non-zero} \}} 
\end{equation}
is finite, and more precisely there exist $c>0$ and $A\geq 0$ such
that for $\ell>\ell(d)$, the set
\begin{equation}\label{eq-field-def}
\bigcup_{[k_1:k]=d}{\{t\in U(k_1)\,\mid\,\text{ there exists $e\not=0$
    in $\mathcal{A}_t[\ell]$ with } [k_1(e):k_1]\leq
  c\ell/(\log\ell)^A\}}
\end{equation}
is finite, where $k_1(e)$ is the field generated by coordinates of
$e$.
\end{thm}

\begin{rem}
  The ``Strong Uniform Boundedness Conjecture'' (which is a theorem
  due to Merel~\cite{merel} in the case $g=1$) makes the much stronger
  prediction that the set
 $$
\bigcup_{[k_1:k]=d}{\{t\in U(k_1)\,\mid\,
  \text{$\mathcal{A}_t[\ell](k_1)$ is non-zero} \}} 
$$
is \emph{empty} for $\ell$ large enough (depending on $d$), and indeed
that this holds even when $U$ is replaced by the entire moduli space
of abelian $g$-folds!
\end{rem}

% In both special cases above, the case $d=1$ only requires the
%verification of condition (i), i.e. the growth in genus of a family of
%covers of $U$.  The genus is a purely toplogical, or even
%combinatorial invariant of the cover, and can be bounded below by
%group-theoretic means.  By contrast, the results of the present paper
%demand that the {\em gonality} of the covers $U_i$ be bounded below;
%the gonality requires

\subsection{Outline of the proofs and of the paper}

We now briefly summarize the basic ideas in the proofs.  For
Theorem~\ref{th:main}, the argument is quite short but involves a
rather disparate combination of ideas. We proceed in four steps, whose
combination is rather surprising. First, using the esperantist
property and a result of Kelner~\cite{kelner}, we show that the genus
of the smooth projective models $C_i$ of the $U_i$ goes to infinity;
second, we invoke comparison principle between the first eigenvalue of
the Cayley-Schreier graphs attached to the covers $U_i$ and the first
Laplace eigenvalue on the Riemann surface $U_{i,\C}$ (this goes back
to Brooks~\cite{brooks2, brooks} and Burger~\cite{burger}); next, we
combine these facts to infer, by means of a theorem of Li and
Yau~\cite{li-yau}, that the \emph{gonality} of $U_i$ tends to
infinity; finally, a result of Abramovich and Voloch~\cite{av} or
Frey~\cite{frey} (which involve Faltings' Theorem~\cite{faltings2} on
rational points on subvarieties of abelian varieties) gives the
desired uniformity for points of bounded degree.
\par
For Theorems~\ref{th:sp} and \ref{th:sl2sl2}, it is easy to describe a
suitable family of covers for which the conclusion of
Theorem~\ref{th:main} leads to the desired conclusion: they are
constructed from congruence quotients of the image $\Gamma$ of the
relevant monodromy representation. The main difficulty is to prove
that these covers satisfy the esperantist property. In
Section~\ref{sec-sources}, we explain different results which provide
this property. The easiest case is when $\Gamma$ has Property (T) of
Kazhdan (see~\cite{bhv}), which turns out to happen in the special
case of Corollary~\ref{cor:concrete} when $g\geq 2$ (because of a
result of J-K Yu~\cite{yu}). However, most of the time this property
is either not known, or not true. When $\Gamma$ is of infinite index
in a lattice in its Zariski-closure, which is a crucial case, we
derive the esperantist property by an appeal to the remarkable recent
results concerning expansion in linear groups over finite
fields. These begin with Helfgott's breakthrough treatment of $\SL_2$
(see~\cite{helfgott}), and further cases are due to Gill and
Helfgott~\cite{h-g}, Breuillard-Green-Tao~\cite{bgt} and particularly
Pyber-Szab\'o~\cite{psz}. In the case of Theorem~\ref{th:sl2sl2},
there are very simple concrete examples (due to Nori~\cite{nori})
where the image of the relevant representation $\rho$ is not a lattice
(see Example~\ref{ex:nori}), so that these new developments are
absolutely essential.
\par
In the case of Theorem~\ref{th:torsion}, we also require a general
result that states, roughly speaking, that the image of the monodromy
representation modulo $\ell$ are ``almost'' perfect groups generated
by elements of order $\ell$, which we prove in
Section~\ref{sec:general}. The esperantist property in such a case is
also obtained from the work of Pyber and Szab\'o.
\par
We think that this paper raises a number of interesting questions. In
the final Section, we raise a few of them, and also make a few
additional comments.  Finally, in two appendices, we record some
necessary facts about comparison between combinatorial and analytic
Laplacians (Appendix A) and semisimple approximation of linear groups
over finite fields (Appendix B) which are either difficult to find in
the form we need, or not fully spelled out, in the published
literature.
\par
\bigskip
\par
\textbf{Notation.}  As usual, $|X|$ denotes the cardinality of a set,
or, if $X$ is a graph, the cardinality of its vertex set.
\par
By $f\ll g$ for $x\in X$, or $f=O(g)$ for $x\in X$, where $X$ is an
arbitrary set on which $f$ is defined, we mean synonymously that there
exists a constant $C\geq 0$ such that $|f(x)|\leq Cg(x)$ for all $x\in
X$. The ``implied constant'' refers to any value of $C$ for which this
holds. It may depend on the set $X$, which is usually specified
explicitly, or clearly determined by the context.
\par
If $U/k$ is an algebraic curve over a number field, we denote by
$U_{\C}$, or sometimes $U(\C)$, the associated Riemann surface, with
its complex topology.  If $k$ is a number field, we write $\OO_k$ for
its ring of integers, and 
$$
\bigcup_{[k_1:k]=d}{(\cdots)}
$$ 
denotes a union over all extensions $k_1/k$ of degree $d$. When
considering \'etale or topological fundamental groups, we often omit
explicit mention of a basepoint.
\par
\medskip
\par
\textbf{Acknowledgments.}  We warmly thank J. Achter for sending us
the paper~\cite{masser} of D. Masser, which motivated us to extend the
results of~\cite{eehk} to points of bounded degree. Thanks also to
M. Burger for discussions and clarifications concerning the links
between expanders and Laplace eigenvalues, and to A. Gamburd for
useful discussions at an earlier stage of this project. We also thank
P. Sarnak for pointing out to us the paper~\cite{z} of P. Zograf and
for other enlightening remarks, and J. Bourgain for pointing out that
the results of E. Hrushovski~\cite{hru} are also applicable to obtain
gonality growth.
\par
The first-named author's work was partially supported by NSF-CAREER
Grant DMS-0448750 and a Sloan Research Fellowship.

\section{Growth of gonality in expanding families of coverings}

Theorem~\ref{th:main} is obtained by combining known finiteness
statements from arithmetic geometry and a ``geometric growth'' theorem
of independent interest, which we state first.  For a smooth curve
$X/k$, $k$ a number field, We recall that the \emph{gonality}
$\gamma(X)$ is the minimal degree of a dominant morphism from $X_{\C}$
to $\P^1_{\C}$.

\begin{thm}[Growth of genus and gonality]
\label{th:gonalitymain}
Let $U/\C$ be a smooth connected algebraic curve over $\C$. Let
$(U_i)_{i\in I}$ be an infinite family of \'etale covers of $U$, and
let $C_{i}$ be the smooth projective model of $U_i$. Let $S$ be a
fixed finite symmetric generating set of the topological fundamental
group $\pi_1(U,x_0)$ for some fixed $x_0\in U$. Assume that the family
of Cayley-Schreier graphs $C(N_i,S)$ associated to the finite quotient
sets
$$
N_i=\pi_1(U,x_0)/\pi_1(U_{i},x_i),\quad\quad x_i\in
U_{i}\text{ some point over } x_0,
$$
is an esperantist family with constants $c>0$ and $A\geq 0$.
\par
\emph{(a)} The genus $g(C_i)$ tends to infinity, i.e., for any $g\geq
0$, there are only finitely many $i$ for which the genus $g(C_i)$ of
$C_i$ is $\leq g$.
% More precisely, if $\lambda>0$ is a
% lower bound for $\lambda_1(C(N_i,S))$, we have
% $$
% \liminf_{i\ra +\infty}{\frac{g(C_i)}{\lambda |N_i|}}>0.
% $$
\par
\emph{(b)} The gonality $\gamma(U_i)$ also tends to infinity. In fact,
there exists a constant $c'>0$ such that
\begin{equation}\label{eq-quant-gonality}
\gamma(U_i)\geq c'|N_i|/(\log 2|N_i|)^{2A}
\end{equation}
for all $i$.
% We have
% $$
% \lim_{i\ra +\infty}{\gamma(U_i)}=+\infty,
% $$
% where $\gamma(U_i)$ denotes the gonality of $U_i$.
% and more precisely
% $$
% \liminf_{i\ra +\infty}{\frac{\gamma(C_i)}{\lambda g(C_i)}}\geq 4\pi.
% $$
% r any $\gamma\geq 1$, there are only finitely many $i$
% for which the gonality $\gamma(U_i)$ is $\leq \gamma$.
\end{thm}

Note that (a) is certainly a necessary condition for (b) to be
true. However, for the proof of (b), we will need a form of (a), so we
have included the statement separately. One may already see that, by
applying the Mordell Conjecture (proved by Faltings~\cite{faltings}),
one deduces from (a) the case $d=1$ of Theorem~\ref{th:main}. In the
setting of Theorem~\ref{th:sp}, this argument provides a new proof of
some of our results in~\cite{eehk}.
\par
In order to deduce the full force of Theorem~\ref{th:main} from
Theorem~\ref{th:gonalitymain}, we use a well-known and very deep
result of Faltings and Frey, which roughly states that curves over
number fields only have infinitely many points of bounded degree for
``obvious'' reasons:
% first step, accordingly, is to
% prove this.  Note also that the case $d=1$ of Theorem~\ref{th:main}
% follows already from (a), by a direct application of Faltings's proof
% of the Mordell conjecture~\cite{faltings}. 

% Our finiteness theorems ultimately derive from arithmetic geometry, in
% particular from a consequence of Faltings' proof of one of Lang's
% conjectures on rational points of higher-dimensional algebraic
% varieties.   For
% instance, for any hyperelliptic curve
% $$
% y^2=f(x),
% $$
% with $f$ squarefree of degree $\geq 3$, the gonality is $2$. Note that
% gonality, defined in this manner, is a birational invariant. By
% contrast with the genus, gonality is \emph{not} a purely topological
% invariant; holomorphic curves of the same genus, which form a single
% homeomorphism class, can have gonality ranging from $2$ to $\lfloor
% \frac{g+3}{2} \rfloor$.

\begin{thm}[Faltings, Frey]\label{th:criterion}
  Let $k$ be a number field, and let $X/k$ be a smooth geometrically
  connected algebraic curve. For any positive integer $d$ such that
  $\gamma(X)>2d$, the set
$$
\bigcup_{[k_1:k]=d}{X(k_1)}
$$
is finite, i.e., there are only finitely many points of $X$ defined
over an extension of $k$ of degree at most $d$.
\end{thm}

This is Proposition 2 in~\cite{frey}: Frey shows that the existence of
infinitely many points over extensions of $k$ of degree $d$ implies
the existence of a non-trivial $k$-rational map of degree $\leq 2d$ to
$\P^1_k$, using the main theorem of Faltings on rational points of
abelian varieties~\cite{faltings}. (Alternatively, it was observed by
Abramovich and Voloch~\cite{av} that this follows from a result of
Abramovich and Harris~\cite[Lemma 1]{ah} and the theorem of Faltings.)
\par
We now begin the proof of Theorem~\ref{th:gonalitymain}. First, we
obtain (a) quite quickly. The basic intuition here is that a graph
embedded in a surface of bounded genus can not have a large Cheeger
constant. Precisely:

\begin{lem}\label{lm:genus-growth}
Under the condition of Theorem~\ref{th:gonalitymain}, we have
$$
g(C_i)\gg \frac{|N_i|}{(\log 2|N_i|)^A}
$$
for all $i$, where $A$ is the constant appearing in
Definition~\ref{def:esperantist}. 
\end{lem}

\begin{proof}
  Following the ideas of Brooks~\cite[\S 1]{brooks2} and of
  Burger~\cite{burger} (see also~\cite[p. 50]{lubotzky}, and
  Appendix~A), it is known that there is a suitable symmetric set of
  generators $S_0$ of $\pi_1(U,x_0)$ such that the Cayley-Schreier
  graph $\Gamma_i=C(N_i,S_0)$ with respect to $S_0$ may be
  \emph{embedded} in the Riemann surface $U_i(\C)$, hence in
  $C_i(\C)$, for all $i\in I$. Moreover, this family $(\Gamma_i)$ is
  still an esperantist family (see for
  instance~\cite[Th. 4.3.2]{lubotzky} for this standard fact).
\par
Now, a beautiful result of Kelner~\cite[Th. 2.3]{kelner} shows that
the first non-zero eigenvalue $\lambda_1(\Gamma_i)$ satisfies
\begin{equation}\label{eq:kelner}
  \lambda_1(\Gamma_i)\ll \frac{\max(g(C_i),1)}{|\Gamma_i|}
\end{equation}
where the implied constant depends only on the degree of the graph,
and the result follows from Definition~\ref{def:esperantist}.
\end{proof}

Using this, we now go to the second part of
Theorem~\ref{th:gonalitymain}, the lower bound on the gonality.

\begin{proof}[Proof of \emph{(b)}]
  In a slight abuse of notation, we use $U_i$ (resp. $C_i$) here to
  refer to the Riemann surfaces $U_i(\C)$, $C_i(\C)$. Following ideas
  of Zograf~\cite{z} and Abramovich~\cite[\S 1]{a}, we first use a
  result of Li and Yau to connect the gonality of $U_i$ to its genus
  its first Laplace eigenvalue.
\par
By (a), the universal cover of the (possibly open) curve $U_i$ is the
hyperbolic plane $\mathbb{H}$ for all but finitely many $i$, and we
can exclude these exceptions from consideration. We can therefore
represent $U_i$ as a quotient
$$
U_i\simeq G_i\backslash \mathbb{H},
$$
where $G_i$ is a discrete subgroup of $\mathrm{PSL}_2(\R)$. From
$\mathbb{H}$, the Riemann surface $U_i$ also inherits the hyperbolic
metric and its associated area-element $d\mu=y^{-2}dxdy$. The
hyperbolic area of $U_i$ is finite (since $U_i$ differs from the
compact curve $C_i$ by finitely many points, i.e., it is a Riemann
surface of finite type). Moreover, the Poincar\'e metric also induces
the Laplace operator $\Delta$ on the space $L^2(U_i,d\mu)$. Thus one
can define the invariant
\begin{align}
\lambda_1(U_i)&=
\inf\Bigl\{
\frac{\langle \Delta \varphi,\varphi\rangle}{\|\varphi\|^2}
\,\mid\, \varphi\text{ smooth and }
\int_{U_i}{\varphi(x)d\mu(x)}=0
\Bigr\}
\nonumber\\
&=\inf\Bigl\{
\frac{\displaystyle{\int_{U_i}{\|\nabla \varphi\|^2d\mu}}}{\|\varphi\|^2}
\,\mid\, \varphi\text{ smooth and }
\int_{U_i}{\varphi(x)d\mu(x)}=0
\Bigr\},
\label{eq:lambdaone}
\end{align}
where 
$$
\nabla \varphi=y^2(\partial_x\varphi,\partial_y\varphi)\,:\,
U_i\ra \C^2
$$ 
is the gradient of $\varphi$, computed with respect to the hyperbolic
metric. It is known that $\lambda_1(U_i)$ is either equal to $1/4$ or
to the first non-zero eigenvalue of the laplacian $\Delta$ acting on
$L^2$-functions on $U_i$. 
\par
It follows from the results of Li and Yau that
\begin{equation}\label{eq:liyau}
\gamma(U_i)\geq \frac{1}{8\pi}
\lambda_1(U_i)\mu(U_i).%\int_{U_i}{d\mu}.
\end{equation}
\par
More precisely, writing $V_c(2,U_i)$ for the conformal area of $U_i$
(as defined in~\cite[\S 1]{li-yau}) the easy bound~\cite[Fact 1, Fact
2]{li-yau}
$$
V_c(2,U_i)\leq \gamma(U_i)V_c(2,\mathbb{S}^2)=4\pi \gamma(U_i)
$$
and the key inequality \cite[Th. 1]{li-yau}
$$
\lambda_1(U_i)\mu(U_i)\leq 2V_c(2,U_i)
$$
combine to give (\ref{eq:liyau}).  Note that, although Li and Yau
assume that the surfaces involved are compact, Abramovich~\cite{a}
explains how the inequality~(\ref{eq:liyau}) extends immediately to
finite area hyperbolic surfaces.
\par
By the Gauss-Bonnet theorem (for finite-area hyperbolic surfaces, see,
e.g.,~\cite[Th. B]{rosenberg}) we have
$$
\mu(U_i)=-2\pi \chi(U_i)=-2\pi(\chi(C_i)-|C_i-U_i|)\geq 4\pi(g(C_i)-1)
$$
where $\chi(\cdot)$ denotes the Euler-Poincar\'e characteristic,
hence~(\ref{eq:liyau}) leads to
% (see, e.g.,~\cite[Th. B]{rosenberg}); together with the formulas
% $$
% \chi(U_i)=\chi(C_i)-|C_i-U_i|\leq \chi(C_i)=2(1-g(C_i)),
% $$
\begin{equation}\label{eq:gonality-lower}
\gamma(U_i)\geq \frac{1}{8\pi}\lambda_1(U_i)\ (-2\pi\chi(C_i))\geq 
2\lambda_1(U_i)(g(C_i)-1).
\end{equation}
\par
In turn, Lemma~\ref{lm:genus-growth} now gives
$$
\gamma(U_i)\gg \lambda_1(U_i)\frac{|N_i|}{(\log 2|N_i|)^A}
$$
using the esperantist condition (for all but finitely many $i$). 
\par
Finally, the comparison principle of Brooks~\cite{brooks} and
Burger~\cite{burger,burger2} relates $\lambda_1(U_i)$ to the
combinatorial laplacian of the Cayley-Schreier graphs: by~\cite[\S 3,
Cor. 1]{burger2}, there exists a constant $c>0$, depending only on $U$
and on $S$, such that
\begin{equation}\label{eq-burger}
\lambda_1(U_i)\geq c \lambda_1(\Gamma_i)
\end{equation}
for all $i$ (Brooks and Burger state their result for \emph{compact}
Riemannian manifolds, but they also both mention that they remain
valid for finite-area Riemann surfaces; in Theorem~\ref{th:burger} in
Appendix~A, we sketch the extension using Burger's method; see also
the recent extension to include infinite-covolume situations by
Bourgain, Gamburd and Sarnak~\cite[Th. 1.2]{bgs}, although the latter
does not state precisely this inequality). Using once more the
esperantist property, we obtain
$$
\gamma(U_i)\gg \frac{|N_i|}{(\log 2|N_i|)^{2A}}
$$
for all $i$ (we can adjust the implied constant to make the inequality
valid for any finite exceptional set), which is the
conclusion~(\ref{eq-quant-gonality}). 
% , which leads to the desired conclusion (even in the
% stronger form that $\gamma(U_i)\gg g(U_i)$ for all but finitely many
% $i$).
\end{proof}

\begin{rem}  
  Zograf~\cite[Th. 5]{z} first showed the relevance of arguments of
  differential geometry of Yang and Yau~\cite{yy} to prove gonality
  bounds for modular curves.  The result of Li and Yau is similar to
  that of~\cite{yy} and both are remarkable in that they prove a lower
  bound for the degree of any conformal map $C_i\ra \mathbb{S}^2$, in
  terms of the hyperbolic area of $C_i$ and the first Laplace
  eigenvalue. These arguments are highly ingenious, involving an
  application of a topological fixed-point theorem to find a suitable
  test function in order to estimate $\lambda_1$.
  Abramovich~\cite{a}, independently, also applied~\cite{li-yau} to
  modular curves.
\end{rem}

\begin{rem}
  There is an intriguing similarity between the proof of
  Theorem~\ref{th:gonalitymain} and a beautiful recent result of
  Gromov and Guth~\cite[Th. 4.1, 4.2]{gromov-guth}. Roughly speaking,
  their result implies that for any family $(M_i\ra M)$ of Galois
  coverings of hyperbolic $3$-manifolds, and for any knot $K_i$ such
  that $K_i$ is the ramification locus of a map $M_i\ra \mathbb{S}^3$
  of degree $3$ (which exists by results of Hilden and Montesinos in
  $3$-dimensional topology), the \emph{distortion} of $K_i$
  (see~\cite[\S 4]{gromov-guth} for the definition) goes to infinity
  if the Cayley graphs of $\pi_1(M)/\pi_1(M_i)$ (with respect to a
  fixed symmetric generating set of $\pi_1(M)$) form an expander
  family. However, their quantitative lower bound shows that it is
  enough to consider an \emph{esperantist} family.
\end{rem}

\section{Sources of expansion}\label{sec-sources}

We explain in this section some results which ensure that various
families of Cayley-Schreier graphs are esperantist families. Some are
quite classical, while others are very recent developments. We try to
present these in an understandable manner for non-specialists.
\par
\medskip
\par
-- If $G$ is a finite-index subgroup in $\barG(\Q)\cap \GL_m(\Z)$,
where $\barG\inj \GL_m$ is a semisimple algebraic subgroup, defined
over $\Q$, and $\barG$ has real rank at least $2$ (examples include
$\GL_n$, $n\geq 3$, or $\Sp_{2g}$, $g\geq 2$, so that $G$ can be a
finite-index subgroup of $\SL_n(\Z)$, $n\geq 3$, or of $\Sp_{2g}(\Z)$,
$g\geq 2$) and $S$ is an arbitrary finite set of generators of $G$,
then the family of Cayley graphs $C(\Gamma_i,S)$ of \emph{all} finite
quotients $\Gamma_i$ of $G$ is an expander family, and hence an
esperantist family. This is because $G$ has Property (T) of Kazhdan
(see~\cite{bhv} for a full treatment of Property (T)).
\par
\medskip
\par
-- Our families of graphs will most often be based of the following
type: we have a finitely-generated discrete group $G\subset
\barG(\Q)\cap \GL_m(\Z)$ which is \emph{Zariski-dense} subgroup of a
semisimple algebraic group $\barG\inj \GL_m$ (defined over $\Q$), and
we consider the Cayley graphs, with respect to a fixed symmetric set
of generators $S$, of the congruence quotients
$$
G_{\ell}=G/\ker(G\ra \barG(\F_{\ell}))\subset \barG(\F_{\ell})
$$
(which are well-defined for almost all $\ell$ after fixing an integral
model of $\barG$).  Thus $G$ might be of infinite index in the lattice
$\barG(\Q)\cap \GL_m(\Z)$ (the ``thin case'', in the terminology
suggested by~\cite{bgt}), and neither Property (T) nor automorphic
methods are applicable. There has however been much recent progress in
understanding expansion properties of Cayley graphs in this
situation. We will use mostly the following very general criterion of
Pyber and Szab\'o~\cite{psz}:

\begin{thm}[Pyber-Szab\'o]\label{th:psz}
  Let $m\geq 1$ be fixed, let $(G_{\ell})$ be a family of subgroups of
  $\GL_m(\F_{\ell})$ indexed by all but finitely many prime numbers,
  and let $S_{\ell}$ be symmetric generating sets of $G_{\ell}$ with
  bounded order, i.e. $|S_{\ell}|\leq s$ for all $\ell$ and some
  $s\geq 1$. Then, if the groups $G_{\ell}$ are all non-trivial
  perfect groups and are generated by their elements of order $\ell$,
  the family of Cayley graphs $(C(G_{\ell},S_{\ell}))_{\ell}$ is an
  esperantist family.
\end{thm}

This follows from~\cite[Th. 8]{psz}: Pyber and Szab\'o show, under
these assumptions, that the diameters $d_{\ell}$ of the Cayley graphs
$\Gamma_{\ell}=C(G_{\ell},S_{\ell})$ satisfy
\begin{equation}\label{eq-diameter}
d_{\ell}\ll (\log |G_{\ell}|)^{M(m)}
\end{equation}
where $M(m)\geq 0$ depends only on $m$, and we can then apply a
general bound of Diaconis and Saloff-Coste~\cite[Cor. 1]{dsc}: for any
finite group $G$ with symmetric generating set $S$, we have
\begin{equation}\label{eq-dsc}
\lambda_1(C(G,S))\geq \frac{1}{|S|\mathrm{diam}(C(G,S))^2},
\end{equation}
which translates here to
$$
\lambda_1(\Gamma_\ell)\gg \frac{1}{(\log |G_{\ell}|)^{2A}},
$$
proving the esperantist property. 
\par
The flexibility and generality of the Pyber-Szab\'o theorem will be
important in Theorem~\ref{th:torsion}. We should however mention some
earlier results of similar type. The first breakthrough was Helfgott's
proof of an estimate like~(\ref{eq-diameter}) for
$G_{\ell}=\SL_2(\F_{\ell})$ (see~\cite[Main th.]{helfgott}
and~\cite[Cor. 6.1]{helfgott}). This was generalized to $\SL_3$ by
Helfgott, to $\SL_n$ (with some restriction) by Gill and
Helfgott~\cite{h-g}, and to all cases where
$G_{\ell}=\barG(\F_{\ell})$ for a simple split Chevalley group
$\barG/\Z$ by Pyber-Szab\'o~\cite{psz} and (independently)
Breuillard-Green-Tao~\cite{bgt}. This last result could be used,
instead of Theorem~\ref{th:psz}, for all applications except
Theorem~\ref{th:torsion}. In Section~\ref{sec-beyond}, we will make
some further comments concerning recent extensions of these results
which imply that many families of Cayley graphs are in fact expander
families.
\par
\medskip
\par
-- The following simple observation will be used to pass from Cayley
graphs to Cayley-Schreier graphs:

\begin{prop}[Big quotients remain esperantist]\label{pr-big-quotients}
Let $(C(G_i,S_i))$ be an esperantist family of Cayley graphs, and let
$N_i\subset G_i$ be subgroups of $G_i$ such that
$$
\log (2[G_i:N_i])\geq \delta \log (2|G_i|)
$$
for some $\delta>0$ independent of $i$. Then the family of
Cayley-Schreier graphs $(C(G_i/N_i,S_i))$ is also an esperantist
family.
\end{prop}

\begin{proof}
  Because $\Gamma_i=C(G_i/N_i,S_i)$ is a quotient of the Cayley graph
  $\tilde{\Gamma}_i=C(G_i,S_i)$, the Laplace eigenvalues of
  $\tilde{\Gamma}_i$ appear among those of $\Gamma_i$. Thus
$$
\lambda_1(\Gamma_i)\geq \lambda_1(\tilde{\Gamma}_i)
\geq \frac{c}{(\log 2|G_i|)^{A}}\geq 
\frac{c\delta^A}{(\log 2|G_i/N_i|)^A},
$$
which proves the esperantist property of the $(\Gamma_i)$. 
\end{proof}

% are bounded by $(\log |G_{\ell}|)^{M(m)}$ for some constant $M(m)\geq
% 0$ independent of $\ell$, and apply the bounds of .  group There are two justifications
% for Definition~\ref{def:esperantist}. The first is that (by work of
% Diaconis and Saloff-Coste~\cite[Cor. 1]{dsc}), any family of Cayley
% graphs of bounded degree satisfying~(\ref{eq:limit-size}) and having
% poly-logarithmic diameter
% $$
% \mathrm{diam}(\Gamma_i)\ll (\log |\Gamma_i|)^{A}
% $$
% for some $A\geq 0$ satisfies
% $$
% \lambda_1(\Gamma_i)\gg \frac{1}{(\log |\Gamma_i|)^{2A}},
% $$
% and hence is an esperantist family. 
% \par
% The other reason is that, in the recent
% works~\cite{helfgott,bg,bgt,psz} proving that Zariski-dense subgroups
% of arithmetic groups have Property $(\tau)$ with respect to congruence
% subgroups, the current methods rely on first proving a generalization
% of the growth theorem of Helfgott~\cite[Main th.]{helfgott}, which
% \emph{by itself} implies easily that the corresponding families of
% graphs are esperantist (see~\cite[Cor. 6.1]{helfgott}). To prove the
% stronger statement that these graphs are expanders typically involves
% significant further ideas (as in~\cite{bg}), and it is therefore of
% interest to know that this strengthening is not required for
% our applications. 
% \par
% In fact, in most cases in this paper, the expansion has not yet been
% proved in the published literature.  Our main tool will then
% be the following theorem of Pyber and Szab\'o:
% This is an immediate consequence of 
\par
\medskip
\par
We are now ready to give the proofs of the arithmetic applications of
Theorem~\ref{th:main} described in the introduction, referring to
Section~\ref{sec-beyond} for more discussion and comparison of the
esperantist condition with the more common expander condition. But
before, here are two direct applications of
Theorem~\ref{th:gonalitymain}, which seem enlightening.

\begin{prop}[Existence of towers with increasing gonality]
  Let $X_0/\C$ be a compact connected Riemann surface of genus $g\geq
  2$. There exists a tower
$$
\cdots\ra X_{n+1}\ra X_n\ra X_{n-1}\ra \cdots \ra X_0
$$
of \'etale Galois coverings such that the gonality of $X_n$ tends to
infinity as $n\ra +\infty$.
\end{prop}

\begin{proof}
  By Theorem~\ref{th:gonalitymain}, it is enough to construct a tower
  $(X_n)$ of this type in such a way that the Cayley graphs of the
  Galois groups form an expander with respect to a fixed set of
  generators of $\pi_1(X_0)$. This is possible because $\pi_1(X_0)$ is
  sufficiently big (since $g\geq 2$) to have a quotient which is a
  discrete group with Property $(T)$, e.g., $\SL_3(\Z)$
  (see~\cite[Cor. 6]{brooks} and the first item in this section).
\end{proof}

The second example illustrates that sometimes esperantism is the best
that can be hoped for:

\begin{prop}
  Let $U=\P^1-\{0,1,\infty\}$, let $x_0\in U$ be any point, and let
  $S=\{a_0^{\pm 1},a_1^{\pm 1}\}$ be the generating set of
  $\pi_1(U,x_0)$ where $a_i$ is a loop around $i$. There exists a
  family of \'etale Galois covers $(U_{m,k}\ra U)_{m\geq 3, k\geq 2}$
  with Galois group $\Gal(U_{m,k}/U)\simeq \SL_m(\Z/k\Z)$ such that
  the Cayley graphs $C(\Gal(U_{m,k}/U),S)$ form an esperantist, but
  not an expander, family. In particular, for all $N\geq 1$, only
  finitely many of $U_{m,k}$ have gonality $\leq N$.
\end{prop}

\begin{proof}
  Kassabov and Riley~\cite[Th. 1.1]{kr} give explicit two-element
  generating sets $S_m=\{a_m,b_m\}$ of $\SL_m(\Z)$ for $m\geq 3$ such
  that the diameter of $\SL_m(\Z/k\Z)$ with respect to $S_m\cup
  S_m^{-1}$ satisfies
$$
d_{m,k}\leq 3600\log (|\SL_m(\Z/k\Z)|),
$$
so that the family of Cayley graphs $(C(\SL_m(\Z/k\Z),S_m\cup
S_m^{-1}))_{m\geq 3,k\geq 2}$ is an esperantist family, but the
Cheeger constant is $\ll 1/m$ (an observation of Y. Luz), so that it
can not be an expander family.
\par
Since $\pi_1(U,x_0)$ is a free group generated by $a_0$ and $a_1$, we
have surjective homomorphisms
$$
\phi_{m,k}\ 
\begin{cases}
  \pi_1(U,x_0)\ra \SL_m(\Z/k\Z)\\
  a_0\mapsto a_m\mods{k}\\
  a_1\mapsto b_m\mods{k},
\end{cases}
$$
and we can define $U_{m,k}$ as the covering of $U$ associated to the
kernel $\ker(\phi_{m,k})$. The Cayley graph $C(\Gal(U_{m,k}/U),S)$ is
isomorphic to $C(\SL_m(\Z/k\Z),S_m)$, hence the result follows.
\end{proof}

\section{First arithmetic applications}\label{sec:application}

We give here the proofs of Theorems~\ref{th:sp}
and~\ref{th:sl2sl2}. Theorem~\ref{th:torsion}, which requires more
preparatory work, is considered in the next section.

\subsection{Proof of Theorem~\ref{th:sp} and
  Corollary~\ref{cor:concrete}}

The proof of Theorem~\ref{th:sp} begins like the argument
in~\cite{eehk}.  We denote by $\Gamma \subset \Sp_{2g}(\Z)$ the image
of the topological monodromy homomorphism $\rho$ associated to
$\mathcal{A}$, and we fix once and for all a finite symmetric
generating set $S$ of $\pi_1(U_{\C},x_0)$. 
\par
By hypothesis, $\Gamma$ is Zariski-dense in $\Sp_{2g}(\Z)$, and
therefore suitable forms of the Strong Approximation Theorem
(see~\cite{mvw}) imply that the image $\Gamma_{\ell}$ of the reduction
map
$$
\Gamma\ra\Sp_{2g}(\F_{\ell})
$$
is equal to $\Sp_{2g}(\F_\ell)$ for all but finitely many $\ell$.
\par
The pairs $(\ell,H)$, where $\ell$ is a prime such that
$\Gamma_\ell=\Sp_{2g}(\F_\ell)$ and $H<\Gamma_\ell$ varies in a fixed
set of representatives of the conjugacy classes of maximal proper
subgroup of $\Gamma_{\ell}$ form an infinite countable set $I$.  To
each such pair $i=(\ell,H)$ corresponds an \'etale $k$-covering
$$
U_i \fleche{\pi_i} U
$$
with a graph isomorphism
$$
C(\pi_1(U_{\C},x_0)/\pi_1(U_{i,\C},x_i),S)\simeq C(N_i,S)
$$
where $N_i=N_{\ell,H}=\Gamma_{\ell}/H$ and $x_i$ is some point in
$U_i$ over $x_0$.
\par
In particular, for any finite extension $k_1/k$, we have
$$
\{t\in U(k_1)\,\mid\,
\mathrm{Im}(\bar\rho_{t,\ell})\not\supset
\Gamma_{\ell}%\cap\Gamma_\ell<\Gamma_\ell
\} \ \subset\
\bigcup_{(\ell,H)\in I}{\pi_{\ell,H}(U_{\ell,H}(k_1))}.
$$
\par
Since the set of $H$ for a given $\ell$ is finite, it follows that
Theorem~\ref{th:main} leads to the desired conclusion once we know
that the family of Cayley-Schreier graphs $(C(N_i,S))$ is an
esperantist family. 
\par
We first check~(\ref{eq:limit-size}) holds, which is easy: if the
degree of the $U_i$ over $U$ were bounded by some $N$, the Galois
group of the splitting field of $U_i/U$ would be contained in $S_N$
for all $i$, which is evidently not the case for $\ell$ large enough.
\par
The esperantist property~(\ref{eq:esperanto}) follows from
Theorem~\ref{th:psz} of Pyber-Szab\'o (or from~\cite{bgt}). Precisely,
it is well-known that $\Sp_{2g}(\F_{\ell})$ is perfect for $\ell\geq
5$, which we may assume, and generated by elements of order $\ell$,
and hence we obtain first
$$
\lambda_1(C(\Gamma_{\ell},S)) \geq \frac{c}{(\log |\Gamma_{\ell}|)^A}
$$
for some $c>0$ and $A\geq 0$. Then, since $N_{\ell,H}=\Gamma_{\ell}/H$
and it is known that the index of any maximal subgroup $H$ of
$\Sp_{2g}(\F_{\ell})$ satisfies
$$
|N_{\ell,H}|=[\Sp_{2g}(\F_{\ell}):H]\geq \tfrac{1}{2}(\ell^g-1)
$$
(e.g., from~\cite[Lemma 4.6]{ls} and Frobenius reciprocity), we can
apply Proposition~\ref{pr-big-quotients} to conclude that the family
$(C(N_{\ell,H},S))_{\ell,H}$ is also an esperantist family. This
concludes the proof of Theorem~\ref{th:sp}.

% In fact, it is very likely that it will be proved
% in the near future that it is an expander family, from the recent
% (independent) results of Helfgott~\cite{helfgott}, Gill and
% Helfgott~\cite{h-g}, Breuillard, Green, and Tao~\cite{bgt} and Pyber
% and Szab\'o~\cite{psz}, combined with the methods of
% Bourgain-Gamburd~\cite{bg}. Indeed, it should then follow that the
% family of Cayley-Schreier graphs
% $$
% (C(\Gamma_{\ell},S))_{\ell}
% $$
% is an expander family. Then, for any $(\ell,H)\in I$, we have a graph
% covering 
% $$
% C(\Gamma_{\ell},S)\ra C(N_{\ell,H},S),
% $$
% and it follows formally that
% $$
% \lambda_1(C(N_{\ell,H},S))\geq \lambda_1(C(\Gamma_{\ell},S))
% $$
% (the pullback to $C(\Gamma_{\ell},S)$ of an eigenfunction of the
% combinatorial Laplace operator on $C(N_{\ell,H},S)$ being another
% eigenfunction with the same eigenvalue). This gives the spectral gap,
% so 
% \par
% However, current literature does not contain a proof of the expansion,
% and we now show how to derive
%%%%%%%%%%%%%%
\par
\medskip
\par
We now come to the proof of Corollary~\ref{cor:concrete}.
By~\cite[Prop. 4]{eehk}, we know that, for all $t \in k$ and all
sufficiently large $\ell$ (in terms of $g$), the surjectivity of the
mod $\ell$ Galois representation
$$
\bar{\rho}_{t,\ell}\,:\, \Gal(\bar{\Q}/k(\zeta_\ell))\ra
\Sp(J(\mathcal{C}_t)[\ell]) \cong \Sp_{2g}(\F_{\ell})
$$
implies that $\End_{\C}(J(\mathcal{C}_t))=\Z$. 
\par
Thus it is enough to prove that Theorem~\ref{th:sp} is applicable to
these families of jacobians. This follows from a theorem of
J-K. Yu~\cite[Th. 7.3 (iii), \S 10]{yu}: the image of the monodromy
representation is not only Zariski-dense in $\Sp_{2g}$, in that case,
it is known precisely to be the principal congruence subgroup
$$
\Gamma=\{x\in \Sp_{2g}(\Z)\,\mid\, x\equiv 1\mods{2}\},
$$
which is of finite index in $\Sp_{2g}(\Z)$ (Yu derives this from the
explicit form of the monodromy around each missing point $t\in \C$
with $f(t)=0$; indeed, these are $2g$ transvections, and Yu is able to
compute precisely the group they generate in $\Sp_{2g}(\Z)$). In other
words, for this particular case, we can appeal to Property (T) (as
explained in the beginning of Section~\ref{sec-sources}) instead of
using the Pyber-Szab\'o theorem. (There is a fairly direct and
elementary proof of Property (T) for these groups, due to
Neuhauser~\cite{neuhauser}, based on methods of Shalom.)

\subsection{Proof of Theorem~\ref{th:sl2sl2}}

The idea is quite similar to the proof of Theorem~\ref{th:sp}.  First
of all, for any prime $\ell$, there exists a cover
$$
U_{\ell}\fleche{\pi_{\ell}} U
$$
defined over $k$ such that $U_{\ell}$ parametrizes pairs $(t,\phi)$
where $t\in U$ and 
$$
\phi\,:\, \EE_{1,t}[\ell] \ra \EE_{2,t}[\ell]
$$
is an isomorphism. It follows that for any finite extension $k_1/k$,
we have
$$
\{t\in U(k_1)\,\mid\,  \EE_{1,t}[\ell]\simeq \EE_{2,t}[\ell]\}
\subset \pi_{\ell}(U_{\ell}(k_1)),
$$
so that the theorem will follow from Theorem~\ref{th:main} once we
establish that the family $(U_{\ell})_{\ell}$ has the desired
expansion property.
\par
First of all, we observe that we may assume both $\EE_1$ and $\EE_2$
are non-isotrivial (indeed, since they are not geometrically
isogenous, at most one can be isotrivial; if $\EE_1$ is isotrivial and
$\EE_2$ is not, then after passing to a finite cover of $U$, we can
assume $\EE_1$ is actually constant, and in that case the curves
$U_\ell$ are isomorphic over $\bar{k}$ to the usual modular curves
$X(\ell)$, whose gonality is already known to go to infinity as $\ell
\ra \infty$, see~\cite[Th. 5]{z} and~\cite{a}).
\par
Now consider the monodromy representation
$$
\rho\,:\, \pi_1(U_{\C},x_0)\ra \SL_2(\Z)\times \SL_2(\Z)
$$ 
associated to the ``split'' family $\EE_1\times \EE_2$ of abelian
surfaces. Let $G\subset \SL_2\times \SL_2$ denote the Zariski
closure of the image of $\rho$. Because $\EE_1$ and $\EE_2$ are both
non-isotrivial, we know that $G$ surjects to $\SL_2$ on each factor,
and because $\EE_1$ and $\EE_2$ are geometrically non-isogenous, it
follows by the Goursat-Kolchin-Ribet lemma that in fact we have
$$
G=\SL_2\times\SL_2.
$$ 
%%\textbf{[Reference!]}
\par
Let $V_{\ell}$ be the curve parameterizing triples
$(t,\phi_1,\phi_2)$, where $\phi_i$ are isomorphisms
$$
\phi_i\,:\, \EE_{1,t}[\ell]\fleche{\sim} (\Z/\ell\Z)^2.
$$
\par
Then $V_{\ell} \ra U$ is a Galois covering whose Galois group is
contained in $\SL_2(\F_\ell) \times \SL_2(\F_\ell)$, this containment
being an identity for all but finitely many $\ell$ by strong
approximation.  We have a map $V_{\ell} \ra U_{\ell}$ given by
$$
(t,\phi_1,\phi_2)=(t,\phi_2^{-1}\phi_1),
$$
which, for almost all $\ell$, expresses $U_\ell$ as the quotient of
$V_\ell$ by the diagonal subgroup $\Delta \subset \SL_2(\F_\ell)
\times \SL_2(\F_\ell)$.
\par
The esperantist property for the family $(V_{\ell})$ follows from
Theorem~\ref{th:psz} (since $G_{\ell}$ is perfect for $\ell\geq 5$ and
generated by its elements of order $\ell$). As in the proof of
Theorem~\ref{th:sp}, it also follows easily, using
Proposition~\ref{pr-big-quotients}, for the quotients $(U_{\ell})$.

\begin{exmp}\label{ex:nori}
  Let $c\in\Q$ be a fixed rational number not equal to $0$ or $1$ (for
  instance $c=2$). Consider first the Legendre family
$$
\mathcal{L}\,:\, y^2=x(x-1)(x-\lambda)
$$ 
over $V=\A^1-\{0,1\}$. It is well-known that $(\pm 1)$ times the image
of the associated monodromy representation
$$
\pi_1(V,\lambda_0)\ra \SL_2(\Z)
$$
is the principal congruence subgroup $\Gamma(2)$ of level
$2$~\cite{nori}. Now fix some rational number $c\notin\{0,1\}$ and
take $\EE_1=\mathcal{L}$ and $\EE_2$ defined by
$$
\EE_{2,\lambda}=\mathcal{L}_{c\lambda},
$$
both restricted to a common base $U/\Q$, where
$U=\A^1-\{0,1,c^{-1}\}$. 
\par
These two families are non-geometrically isogenous, and hence our
theorem applies. Its meaning is that, in a very strong sense, the
torsion fields of $\EE_{1,\lambda}$ and $\EE_{2,\lambda}$ tend to be
independent. For instance, for a given degree $d\geq 1$, we find that
for all $\ell\geq \ell_0(d)$ large enough (depending on $d$) the set
$$
\{\lambda \in\bar{\Q}\,\mid\, 
[\Q(\lambda):\Q]\leq d,\quad \Q(\lambda,\EE_{1,\lambda}[\ell])
=\Q(\lambda,\EE_{2,\lambda}[\ell])\}
$$
is finite. 
\par
Furthermore, in that case, Nori~\cite{nori} has shown that the image
of the monodromy representation
$$
\mathrm{Im}\bigl(\pi_1(U_{\C},x_0)\ra \SL_2(\Z)\times \SL_2(\Z)\bigr)
$$
is \emph{not} of finite index in $\SL_2(\Z)\times \SL_2(\Z)$. This
means that our result can not be obtained, in that special case, using
only expansion properties of quotients of lattices. We expect that
this phenomenon is much more general. %%%ADD???
\end{exmp}

\section{General abelian varieties}\label{sec:general}

In this section, we will prove Theorem~\ref{th:torsion}. However,
before doing so, some preliminaries of independent interest are
required. Essentially, these amount to proving that, for an arbitrary
one-parameter family of abelian varieties $\mathcal{A}\ra U$ over a
number field $k$, the Galois groups of the coverings $U_{\ell}$
associated to the kernel of the composition
\begin{equation}\label{eq:covering}
\pi_1(U_\C,x_0)\ra \GL_{2g}(\Z)\ra \GL_{2g}(\F_{\ell})
\end{equation}
``almost'' satisfy the assumptions of the Pyber-Szab\'o theorem for
all but finitely many $\ell$. (Theorems~\ref{th:sp}
and~\ref{th:sl2sl2} used special cases of this fact, which were
obvious from the underlying assumptions.)
\par
Since the proof of this fact requires quite different arguments of
arithmetic geometry than those of the rest of the paper, we state the
conclusion in a self-contained way and use it to prove
Theorem~\ref{th:torsion} before going into the details.

\begin{prop}\label{pr:blackbox}
  Let $k$ be a number field and let $U/k$ a smooth geometrically
  connected algebraic curve over $k$.  Let $\mathcal{A}\ra U$ be an
  abelian scheme of dimension $g\geq 1$, defined over $k$. Then there
  exists a finite \'etale cover $V\ra U$ such that, if we denote by
$$
\mathcal{A}_V=\mathcal{A}\times_U V\ra V
$$
the base change of $\mathcal{A}$ to $V$, the image of the monodromy
action on $\ell$-torsion
$$
\pi_1(V_\C,x_0)\ra \GL_{2g}(\F_{\ell})
$$
is, for all but finitely many primes $\ell$, a perfect subgroup of
$\GL_{2g}(\F_{\ell})$ generated by elements of order
$\ell$.\footnote{\ Note that it is permitted for this subgroup to be
  trivial.}
\end{prop}

Using this, which will be proved later, we can prove
Theorem~\ref{th:torsion}.

\begin{proof}[Proof of Theorem~\ref{th:torsion}]
  As in the previous results, we denote by $U_{\ell}$ the covering of
  $U_{\ell}$ of $U_\C$ corresponding to the kernel of the composition
$$
\pi_1(U_\C,x_0)\ra \GL_{2g}(\Z)\ra \GL_{2g}(\F_{\ell}).
$$
\par
Applying Proposition~\ref{pr:blackbox}, we find that, possibly after
performing a base-change to a fixed finite covering $V\ra U$, the
image $G^0_{\ell}$ of this representation is, for all but finitely
many $\ell$ (say, for $\ell\geq \ell_0$), a perfect subgroup of
$\GL_{2g}(\F_\ell)$, generated by its elements of order $\ell$ (which
may be trivial). We start by proving~(\ref{eq-exist-torsion}). Since,
clearly, the finiteness of
$$
\bigcup_{[k_1:k]=d}{\{t\in V(k_1)\,\mid\,
  \text{$\mathcal{A}_{V,t}[\ell](k_1)$ is non-zero} \}} 
$$
implies that of 
$$
\bigcup_{[k_1:k]=d}{\{t\in U(k_1)\,\mid\,
  \text{$\mathcal{A}_t[\ell](k_1)$ is non-zero} \}},
$$
we may assume in fact that $V=U$, without loss of generality.
\par
We will now apply Theorem~\ref{th:main} to conclude. Precisely, it is
enough to show that the non-trivial geometrically connected components
of the (possibly disconnected) covers
$$
\mathcal{A}[\ell] \ra U
$$ 
form an esperantist family as $\ell$ varies (where the trivial
connected component is the image of the zero section $0\,:\, U\ra
\mathcal{A}[\ell]$).  
\par
We let $(U_{\ell,i}\ra U_\C)_{\ell,i}$ denote the family of \'etale
covers of $U_\C$ arising as all Riemann surfaces coming from
non-trivial geometrically connected components of $\mathcal{A}[\ell]$
(the index $i$ parametrizes the components for a given $\ell$),
$\ell$ ranging over primes $\geq \ell_0$.
\par
The covering $\mathcal{A}[\ell]_\C \ra U_\C$ corresponds to the
(not-necessarily transitive) action of $\pi_1(U_{\C})$ on
$\mathcal{A}[\ell]$, which factors through the quotient $\pi_1(U_{\C})
\ra G^0_\ell$.  Any component $U_{\ell,i}$ of $\mathcal{A}[\ell]\ra U$
corresponds to an orbit of this action, hence
$$
\pi_1(U_{\ell,i,\C})/\pi_1(U_{\C})\simeq G_{\ell}^0/H_i 
$$
for some subgroup $H_i$ of $G_{\ell}^0$.  Because $G^0_\ell$ is
generated by elements of order $\ell$, it cannot act non-trivially on a
set of size smaller than $\ell$.  Thus, $U_{\ell,i,\C} \ra U_\C$ is
either an isomorphism or 
\begin{equation}\label{eq-deg-ell}
\deg(U_{\ell,i,\C}\ra U_{\C})\geq \ell.
\end{equation}
\par
For those $\ell\geq \ell_0$ such that $G_{\ell}^0$ is non-trivial
(hence of order $\geq \ell$), we can apply the Pyber-Szab\'o Theorem
(Theorem \ref{th:psz}) to deduce that
$$
\lambda_1(C(G^0_\ell,S))\gg \frac{1}{(\log
  |G_{\ell}^0|)^A}\gg \frac{1}{(\log \ell)^A}
$$
for some constant $A$ independent of $\ell$ (the implied constant
depending also on $g$), and then we derive
$$
\lambda_1(C(G_{\ell}^0/H_i,S))\geq \lambda_1(C(G^0_\ell,S))
\gg \frac{1}{(\log \ell)^A}
\geq \frac{1}{(\log |G_{\ell}^0/H_i|)^{A}}
$$
for every proper subgroup $H_i$ of $G_\ell^0$.
% (since $i$ varies in a finite set for each $\ell$).
\par
Then, the remaining $\ell$, as well as the covers $U_{\ell,i,\C}$ for
which $H_i=G_{\ell}^0$, are covers for which $U_{\ell,i}$ is
isomorphic to $U$, and therefore we only need to show that those exist
only for finitely many $\ell$.  Indeed, such geometric components are
parametrized by the group $A(K\C)$, where $K=k(U)$, which might be
infinite (for instance if $A$ is isotrivial, e.g., if it is a product
$B\times U$, where $B/k$ is a fixed abelian variety over $k$). But for
each extension $k_1/k$ of degree $d$, the only geometric components
which can contribute to $\mathcal{A}_t[\ell](k_1)$ are those which are
themselves defined over the compositum $K k_1$.  So what remains is
just to show that
$$
\bigcup_{[k_1:k]=d} A(K k_1)[\ell]=0
$$
for all $\ell$ large enough.  This is immediate by spreading out
$\mathcal{A}$ and $U$ to a model over an open subscheme of $\OO_k$,
and comparing the torsion of $A(K k_1)$ with the torsion of the fiber
of $\mathcal{A}$ over a finite field.
\par
We consider now the finiteness statement~(\ref{eq-field-def}), arguing
simply that if $(t,e)$ is a pair where $t\in U$ and $e$ is a non-zero
torsion point on $\mathcal{A}_t[\ell]$, and if in addition $t\in
U(k_1)$ for some extension $k_1/k$ of degree $d$ and $k_1(e)$ has
degree $d'$ over $k_1$, then $(t,e)$ corresponds to a point $x\in
U_{i,\ell}(k_1(e))$ for one among our auxiliary curves. In particular,
there are only finitely many such pairs for $dd'\leq
2\gamma(U_{i,\ell})$, which by~(\ref{eq-quant-gonality})
and~(\ref{eq-deg-ell}) translates to the finiteness under the
condition
$$
dd'\ll \ell/(\log \ell)^{2A},
$$
as claimed.
\end{proof} 

We now proceed to the proof of Proposition~\ref{pr:blackbox}. We start
with the following general preliminaries from algebraic geometry. If
$K$ is a field, then given an abelian variety $A/K$ and a finite
extension $L/K$, we write $G_\ell$ for the Galois group of
$L(A[\ell])/L$ and $G_\ell^+\leq G_\ell$ for the characteristic
subgroup generated by the $\ell$-Sylow subgroups of $G_\ell$.
\par
The following theorem shows that when $K$ is finitely generated over
$\Q$, there exists an $L/K$ such that $G_\ell^+$ and $G_\ell/G_\ell^+$
are very nicely behaved for almost all $\ell$.  Serre~\cite{serre:vig}
proved it in the special case where $K$ is a number field, and
indicated that the same argument should extend to finitely generated
fields.

\begin{thm}[Semisimple approximation of Galois groups of torsion
  fields]
  Suppose $K/\Q$ is a finitely generated extension and $A/K$ is an
  abelian variety of dimension $g$.  Then there is a finite extension
  $L/K$ and a constant $c=c(K,A)$ depending only on $K$ and $A$ such
  that if $\ell$ is a prime number $\geq c$, then
  % the Galois group $G_\ell$ of the extension $L_\ell=L(A[\ell])$ of
  % $L$ generated by $\ell$-torsion points of $A$, seen as a subgroup
  % of $\Aut(A[\ell])\simeq \GL_{2g}(\F_{\ell})$ to define
  % $G_{\ell}^+$, satisfy:
\begin{enumerate}
\setlength{\itemsep}{0.04in}
\item $G_\ell/G_\ell^+$ has order prime to $\ell$;
  % \item $G_\ell^+$ acts semisimply on $A[\ell]$;
\item there is a semisimple group
  $\barG_\ell\subset\barGL_{2g}/\F_\ell$ such that
  $\barG_\ell(\F_\ell)^+=G_\ell^+$;
  % \item $G_\ell^+$ is perfect and has bounded index in
  %   $S_\ell=G_\ell\cap\barG_\ell(\F_\ell)$;
\item if $S_\ell=G_\ell\cap\barG_\ell(\F_\ell)$, then $G_\ell/S_\ell$ is abelian.
  % \item for some non-trivial connected torus $\barT_\ell/\F_\ell$
  %   and subgroup $C_\ell\leq G_\ell/S_\ell$ there is an embedding
  %   $C_\ell\to\barT_\ell(\F_\ell)$ with bounded cokernel.
\end{enumerate}
\label{th:semisimple}
\end{thm}
\par
The proof we give for the general case was derived from Serre's.  In
particular, the arguments on algebraic subgroups of
$\barGL_{2g}/\F_{\ell}$ and their $\F_\ell$-rational points are
transported essentially unchanged from Serre's paper; the extra
ingredient is that we need to invoke finiteness theorems for \'etale
covers of positive-dimensional varieties over number fields, while
\cite{serre:vig} only needs finiteness theorems for unramified
extensions of the number field itself.

Before embarking on the proof of the theorem we remark on one behavior
of the Galois groups $G_\ell$ under finite base change.  If $L/K$ is
an arbitrary extension, then replacing $L$ by a finite extension
$L'/L$ (e.g.~the extension induced by a finite extension $K'/K$) has
the effect of replacing $G_\ell,G_\ell^+$ by subgroups
$H_\ell,H_\ell^+$ respectively of index at most $[L':L]$.  Since
$G_\ell^+$ has no proper subgroup of index less than $\ell$ (because
it is generated by its $\ell$-Sylow subgroups), we have also
$H_\ell^+=G_\ell^+$ for $\ell>[L':L]$.
\par
Throughout the proof, we will use ``bounded'' as shorthand for
``bounded by a constant which may depend on $K,L,A$ but which is
independent of $\ell$.''

\begin{proof}
  We start by taking $L=K$, but finitely many times throughout the
  proof of the theorem we will replace $L$ by a finite extension
  $L'/L$.  By the remark following the statement of the theorem, as
  far as the groups $G_\ell^+$ are concerned, the effect of
  such a replacement is to increase $c$.  As far as the quotients
  $G_\ell/G_\ell^+$ are concerned, for $\ell\geq c$, they will
  be replaced by subgroups $H_\ell/H_\ell^+\leq G_\ell/G_\ell^+$.
\par
There is a canonical embedding
$G_\ell\to\Aut(A[\ell])\simeq\GL_{2g}(\F_\ell)$, thus we can apply
results of Nori and Serre to the subgroup
$G_\ell^+\leq\GL_{2g}(\F_\ell)$.  As summarized in Appendix~B, we
start by associating to each (finite) subgroup
$G\leq\GL_{2g}(\F_\ell)$ the characteristic subgroup $G^+\leq G$
generated by its unipotent elements, and then we associate to
$G^+\leq\GL_{2g}(\F_\ell)$ an algebraic subgroup
$\barG^+\subseteq\barGL_{2g}$.  One can say quite a bit about
$\barG^+$, especially when $G$ acts semisimply on $A[\ell]$, and as a
result one can also say quite a bit about $G$ and $G^+$.
\par
If $\ell\geq 2g-1$, then Proposition~\ref{prop:order} implies that
$G_\ell/G_\ell^+$ has order prime to $\ell$, thus (1) holds.
\par
For some constant $\ell_1=\ell_1(2g)$, Theorem~\ref{thm:noriB} implies
$G_\ell^+=\barG_\ell^+(\F_\ell)^+$, for $\ell\geq\ell_1$.  If
$\ell\geq \c(K,A)$, then $G_\ell$ acts semisimply on $A[\ell]$
(cf.~Theorem 1 in \cite[VI.3]{fw}), so if $\ell$ also satisfies
$\ell\geq\ell_1$, then Corollary~\ref{cor:noriB} implies
$\barG_\ell=\barG_\ell^+$ is semisimple.  Hence (2) holds.
\par
We suppose for the remainder of the proof that $\ell\geq\ell_1$ and
that $\barG_\ell$ is semisimple, and we write
$\barN_\ell\subseteq\barGL_{2g}$ for the normalizer of $\barG_\ell$.
The fact that $G_\ell$ normalizes $G_\ell^+$ implies it also
normalizes $\barG_\ell$, thus $G_\ell\leq\barN_\ell(\F_\ell)$.  By
Corollary~\ref{cor:tensor}, there is a positive integer $r=r(2g)$
(independent of $\ell$ and $\barG_\ell$) and a faithful representation
$\barGL_m\to\barGL_r$ which identifies the image of $\barG_\ell$ with
the algebraic subgroup of elements in $\barGL_m$ acting trivially on
the subspace of $\barG_\ell$-invariants.  Moreover, the image of
$\barN_\ell$ in $\barGL_n$ stabilizes this space, and for some $s\leq
r$, its action on the space induces a faithful representation
$\barN_\ell/\barG_\ell\to\barGL_s$.
\par
Let $S_\ell=G_\ell\cap\barG_\ell(\F_\ell)$ and let $J_\ell\leq G_\ell$
be a subgroup of minimal index among those such that $S_\ell\leq
J_\ell$ and $J_\ell/S_\ell$ is abelian.  The image of $G_\ell/S_\ell$
in the faithful representation $G_\ell/S_\ell\to \barGL_s(\F_\ell)$
has order prime to $\ell$, thus we can lift it to a faithful
representation $G_\ell/S_\ell\to\GL_s(\C)$ and apply Jordan's Theorem
to infer that $[G_\ell:J_\ell]=[G_\ell/S_\ell:J_\ell/S_\ell]$ is
bounded.  In particular, we will show that the fixed fields
$L(A[\ell])^{J_\ell}$ all lie a single finite extension $L'/L$, so up
to replacing $L$, (3) will hold.
\par  
So far, the argument has paralleled that in \cite{serre:vig} quite
closely.  We now attend to the new features that appear when $K/\Q$
has positive transcendence degree.
\par  
Write $k$ for the largest algebraic extension of $\Q$ contained in
$K$, and write $S$ for $\Spec \OO_k$.  Let $X$ be a smooth scheme
dominant and of finite type over $\Spec(\Z)$ such that the function
field $k(X)$ is $L$ and such that $A$ has good reduction over $L$, and
let $\pi_1^{\et}(X)$ be the \'etale fundamental group of $X$.  For
each prime $\ell$, let $X[1/\ell]$ be the pullback of $X$ to
$\Spec(\Z[1/\ell])$ and let $\pi_1^{\et}(X[1/\ell])$ be the \'etale
fundamental group of $X[1/\ell]$.  The cover $X_\ell\to X[1/\ell]$
induced by the extension $L_\ell/L$ is \'etale because $A$ has good
reduction over $X$, thus $G_\ell$ is a quotient of
$\pi_1^{\et}(X[1/\ell])$.
\par
Let $\pi_1^t(X[1/\ell])$ denote the quotient of
$\pi_1^{\et}(X[1/\ell])$ corresponding to the maximal \'etale cover
$X'\to X[1/\ell]$ which is tamely ramified over $X - X[1/\ell]$.  The
kernel of the quotient map
$\pi_1^{\et}(X[1/\ell])\to\pi_1^t(X[1/\ell])$ is generated by
pro-$\ell$ groups (coming from wild ramification), hence the image in
$G_\ell$ of this kernel lies in $G_\ell^+$ and the quotient
$\pi_1^{\et}(X[1/\ell])\to G_\ell/G_\ell^+$ factors through
$\pi_1^{\et}(X[1/\ell])\to\pi_1^t(X[1/\ell])$.
\par  
Each irreducible component $Z$ of $X - X[1/\ell]$ gives rise to an
inertia group in $G_\ell$; we now show that the images of these
inertia groups in $G_\ell / G_\ell^+$ generate an {\em abelian} group.
\par
If $I\leq G_\ell$ is one such inertia group, then $I^+$ is the unique
$\ell$-Sylow subgroup of $I$ and $I\to I/I^+$ splits, so there is an
embedding $i:I/I^+\to\GL_{2g}(\F_\ell)$ defined up to conjugation by
an element of $I^+$.  If $\ell\geq\c(K)$, there is a connected torus
$\barI'/\F_\ell$ in $\barGL_{2g}$ such that $i(I/I^+)=\barI'(\F_\ell)$
(see~\cite[Section 1.9]{serre:points}).  Moreover, the characters of
the induced representation $\barI'\to\barGL_{2g}$ all have amplitude
at most $2g$ (see the discussion in \cite[Section
2]{hall-transvections}).
\par
One can show that $\barI'\subset\barN_\ell$, and while $\barI'$
depends on our choice of splitting $I/I^+\to I$, the image
$\barI\subset\barN_\ell/\barG_\ell$ of
$\barI'\to\barN_\ell/\barG_\ell$, which we call an inertial torus, is
canonical because $I^+$ lies in $\barG_\ell(\F_\ell)^+$.  Above we saw
that the induced representation $\barI\to\barGL_s$ comes from the
action of $\barI'$ on the subspace of $\barG_\ell$-invariants in the
tensor representation $\barG_\ell\to\barGL_r$, and thus, if
$n=r_1(2g)$ is the constant in the statement of
corollary~\ref{cor:tensor}, then the characters of $\barI'\to\barGL_r$
and $\barI\to\barGL_s$ have amplitude at most $n\cdot 2g$.
\par
The subgroup $\barI'(\F_\ell)\cap J_\ell$ has bounded index in
$\barI'(\F_\ell)$ because $J_\ell$ has bounded index in $G_\ell$, thus
the subgroup of elements in $\barI(\F_\ell)$ which commute with
$J_\ell/G_\ell$ has the same index or smaller.  Therefore, by an
argument involving rigidity of tori (cf.~\cite[\S
2]{hall-transvections}), if $\ell\geq\c(n,2g)$, then $\barI(\F_\ell)$
commutes with $J_\ell/S_\ell$ in $G_\ell/S_\ell$.  In particular,
$J_\ell'=\barI(\F_\ell)J_\ell$ must lie in $J_\ell$ because of how we
chose the latter.  A similar argument shows that any pair of inertial
tori commute, hence the subgroup
$\barT_\ell\subset\barN_\ell/\barG_\ell$ generated by all such tori is
a connected torus and $\barT_\ell(\F_\ell)\leq J_\ell/S_\ell$.
\par
It follows that the image in $G_\ell$ of the kernel of
$\pi_1^{\et}(X[1/\ell])\to\pi_1^{\et}(X)$ lies in $J_\ell$.  In other
words, there is a maximal normal subgroup $J_\ell'\leq G_\ell$
satisfying $S_\ell\leq J_\ell\leq J_\ell'$ such that the quotient
$G_\ell/J_\ell'$ is bounded and $\pi_1^{\et}(X[1/\ell])\to
G_\ell/J_\ell'$ factors through $\pi_1^{\et}(X)$.  But we know (for
instance, by Theorem 2.9 of \cite{hh}) that there are only finitely
many quotients of $\pi_1^{\et}(X)$ of bounded index.  Thus, we can
replace $X$ with some finite \'{e}tale cover $X' \ra X$ (which has the
effect of replacing $L$ with a finite extension of $L'$) and be
assured that $G_\ell/J_\ell'$ is trivial, which is exactly to say that
$G_\ell / S_\ell$ is abelian, as desired.
\end{proof}

We now study the \emph{geometric} Galois group $G^0_\ell = \Gal(L
\bar{k}(A[\ell]) / L \bar{k})$.

%  CHRIS:  original version of this theorem commented out below
%\begin{thm}
%  Suppose $k/\Q$ is a finitely generated field and $K/k$ is a finitely
%  generated regular extension.  If $A/K$ is a non-trivial abelian
%  variety with no isotrivial quotient, then there is a finite extension
%  $L/K$ and a constant $\ell_0(A,K)$ such that $G_\ell^0=G_\ell^+$ and
%  $G_{\ell}^0\not=1$ for all primes $\ell\geq \ell_0$.
%  \label{g0ell}
%\end{thm}

\begin{thm}
  Suppose $k/\Q$ is a finitely generated field and $K/k$ is a finitely
  generated regular extension.  If $A/K$ is an abelian variety of
  dimension $g$, then there is a finite extension $L/K$ and a constant
  $\ell_0(A)$ such that the geometric Galois group $G_\ell^0$ is a
  perfect subgroup of $\GL_{2g}(\F_\ell)$ generated by elements of
  order $\ell$ for all $\ell \geq \ell_0(A)$.
  \label{g0ell}
\end{thm}

This immediately implies Proposition~\ref{pr:blackbox} by taking
$K=k(U)$ the function field of $U$ and $A/K$ the generic fiber of
$\mathcal{A}\ra U$, after noting that the Riemann surface
corresponding to the \'etale cover of $U$ which has function
field $K(A[\ell])$ is the covering  $U_{\ell}$ defined
by~(\ref{eq:covering}). 

\begin{proof}
  As in the previous theorem, we start with a fixed $L/K$ (the
  extension given in the previous theorem), but finitely many times
  throughout the following proof we may replace $L$ with a finite
  extension $L'/L$.  As far as the groups $G_\ell^0$ are concerned,
  the effect of such a replacement is to increase $\ell_0$.  As far as
  the quotients $G_\ell/G_\ell^0$ are concerned, for $\ell\geq\ell_0$,
  they may be replaced by a proper subgroup $H_\ell/H_\ell^0$.
\par
Let $X/k$ be a smooth geometrically-connected variety such that
$K=k(x)$.  After replacing $X$ by an open dense subscheme, we may
suppose that $A$ has good reduction over $X$ (that is, there is an
abelian scheme $\mathcal{A}/X$ whose generic fiber is $A$.) If we
write $\bar{X}=X\times_k\bar{k}$, we have an exact sequence of \'etale
fundamental groups
$$
\pi_1^{\et}(\bar{X})\lto \pi_1^{\et}(X)\lto \Gal(\bar{k}/k)\lto 1.
$$
\par
Up to replacing $k$ by a finite extension, we may suppose $X(k)$ is
non-empty, and thus that this sequence splits.  In particular, if
$N^0\leq\pi_1^{\et}(\bar{X})$ is an open subgroup, then there is an
open subgroup $N\leq\pi_1^{\et}(X)$ such that $N^0=N\cap
\pi_1^{\et}(\bar{X})$ and $N\to\Gal(\bar{k}/k)$ is surjective.  (In
fact, we may just take $N = N^0 \Gal(\bar{k}/k)$ where
$\Gal(\bar{k}/k)$ is viewed as a subgroup of $\pi_1^{\et}(X)$ via the
chosen splitting.)  Moreover, $\pi_1^{\et}(\bar{X})$, being
topologically finitely generated, has only finitely many quotients of
bounded degree, so if we have an infinite family of \'etale covers
$X_\ell\to X$ with bounded geometric monodromy, then we can find a
finite \'etale cover $X'\to X$ such that the pullbacks $X'_\ell\to X'$
all have trivial geometric monodromy.
\par
Theorem 1 of \cite{kl} implies that the intersection
$[G_\ell,G_\ell]\cap G_\ell^0$ has bounded index in $G_\ell^0$ and
Theorem~\ref{th:semisimple} above implies $[G_\ell,G_\ell]\leq
S_\ell$, so $S_\ell\cap G_\ell^0$ also has bounded index in
$G_\ell^0$.  Thus up to replacing $L$ by a finite extension, we can
assume $G_\ell^0\leq S_\ell$.  The index of $G_\ell^+$ in $S_\ell$ is
bounded, thus the index of $G_\ell^0\cap G_\ell^+$ in $G_\ell^0$ is
also bounded, so up to replacing $L$ by a finite extension, we may
suppose $G_\ell^0\leq G_\ell^+$.
\par
Suppose the algebraic envelope $\barG_\ell$ of $G_\ell^+$ is
semisimple and $\barG_\ell(\F_\ell)^+=G_\ell^+$.  We say a subgroup
$G\leq\GL_{2g}(\F_\ell)$ is quasi-simple if its center $Z\leq G$ has
bounded size and if $G/Z$ is a simple group.  If $\ell\geq 5$, then
$G_\ell^+$ is generated by a bounded set $\Sigma_\ell$ of pairwise
commuting quasi-simple subgroups of $G\leq \GL_{2g}(\F_\ell)$ such
that $G^+=G$.  Moreover, for each $G\in\Sigma_\ell$, the index of
$[G,G]$ in $G$ is bounded by $|Z(G)|$, so if $\ell$ is sufficiently
large, then every $G\in\Sigma_\ell$ is perfect.  For every normal
subgroup $N\leq G_\ell^+$, the commutator subgroup $[N,N]$ has bounded
index in $N$.  It is also generated by a subset of $\Sigma_\ell$, so
$[N,N]^+=[N,N]$.  This applies in particular to $G^0_\ell$, so up to
replacing $L$ by a finite extension $L'/L$, we may suppose, for all
$\ell$, that $G^0_\ell$ is perfect and generated by its elements of
order $\ell$.
\end{proof}

\section{Further remarks and questions}\label{sec:conclusion}

We conclude this paper with some remarks and further questions.

\subsection{Various forms of expansion}\label{sec-beyond}

%%%%%%%%%%%%%%%%%%%%
The strongest possible form of expansion for a family $(\Gamma_i)$ of
graphs is given by the expander condition. This is in fact often a
crucial requirement for applications, as seen for instance in the case
of applications of sieve methods in discrete groups with exponential
growth (see~\cite[\S 5.2]{kowalski} for a discussion). Motivated by
these other applications, there has been extensive, and impressive,
work towards a proof that various families of graphs are expanders. In
particular, after many intermediate works, it has been proved by
Salehi-Golsefidy and Varj\'u~\cite{sgv} that the Cayley graphs of
congruence quotients
$$
\Gamma_m=\Gamma/\ker(\Gamma\ra \GL_n(\Z/m\Z))
$$
form an expander family when $m$ runs over squarefree integers and
$\Gamma$ is a finitely-generated subgroup of $\GL_n(\Z)$ which is
Zariski-dense in a semisimple algebraic group $\barG/\Q$ (the case of
$\barG=\SL_2$, for $m$ prime, is due to Bourgain-Gamburd~\cite{bg},
and for $m$ squarefree to Bourgain-Gamburd-Sarnak, while
Varj\'u~\cite{varju} had obtained the result for $\SL_n$).
\par
Note that although many of the graphs we use in this paper are
quotients of graphs of this type, this is not the case for all: in
Theorem~\ref{th:torsion}, we do not have such precise control of the
groups which appear.  This is one reason why we have chosen to
emphasize the weaker esperantist condition for our graphs. Another is
that the starting point in~\cite{sgv} (as in~\cite{bg}
and~\cite{varju}), when $\Gamma$ which is not a lattice (so that
neither Property (T) nor automorphic methods are available), is the
growth theorem of~\cite{psz} or~\cite{bgt}, which by itself implies
immediately the esperantist property. Moreover, in terms of
effectivity, the esperantist property might be much easier to
establish with actual, explicit, constants (we discuss this in
Section~\ref{ssec-effective}).
%%%%%%%%%%%%%%%%%%%%
\par
In the opposite direction, the reader may have noticed that even
weaker expansion conditions than~(\ref{eq:esperanto}) lead to a growth
of gonality (i.e., to Theorem~\ref{th:gonalitymain}): it would be
enough to have a bound of the form
\begin{equation}
  \lambda_1(N_i) \geq \vartheta(i)|N_i|^{-1/2},\quad
\text{with}\quad \lim_{i\ra +\infty}{\vartheta(i)}=+\infty.
\label{beyondesperanto}
\end{equation}
\par
On the other hand, a variant of the Pyber-Szab\'o Theorem proving this
in the context of Theorem~\ref{th:psz} would \emph{not} suffice for
this paper, since we applied the bound to quotients of Cayley graphs,
and the index of subgroups would not in general be large enough to
preserve~(\ref{beyondesperanto}), as in
Proposition~\ref{pr-big-quotients}.
\par
Another condition which would suffice for this paper is
\begin{equation}\label{eq-hru}
\lambda_1(N_i)\gg |N_i|^{-\eps}
\end{equation}
for all $i$ and any $\eps>0$, the implied constant depending on
$\eps$. Indeed, taking $\eps<1/2$ gives~(\ref{beyondesperanto}), and
if $M_i$ is any quotient of $N_i$ with $\log |M_i|\geq c \log |N_i|$
for some fixed constant $c>0$ (as happens in our applications), taking
$\eps<c/4$, say, leads to~(\ref{beyondesperanto}) for $(M_i)$. This
remark may have some applications since, as J. Bourgain pointed out to
us, Hrushovski~\cite[Th. 1.3, Cor. 1.4]{hru} has proved~(\ref{eq-hru})
in many cases related to our applications in this paper.
\par
In another direction, one might wonder whether there is any kind of
structure theorem for ``highly non-expanding'' Cayley-Schreier graphs
with an eigenvalue violating \eqref{beyondesperanto} with
$\vartheta(i)$ being, say, a large constant.  For instance, if such a
graph is associated to a finitely generated group $\Gamma$ acting on a
finite set $S$, what can one say about the composition factors of the
image of $\Gamma$ in $\Aut(S)$?  (Such structure results are known for
covers $(U_i)$ with bounded genus, in the work of
Guralnick~\cite{gura:monodromy} and
Frohardt-Magaard~\cite{froh:frohardtmagaard}; these families, of
course, fail to satisfy~(\ref{eq-diameter}).)

\subsection{Relation with the work of Cadoret and Tamagawa}

Our main result, and its concrete diophantine applications, are
related to recent work of Cadoret and Tamagawa~\cite{ct1,ct2}. Given
an $\ell$-adic representation
$$
\rho\,:\, \pi_1^{\et}(X)\ra \GL_m(\Z_{\ell})
$$
for some ``nice'' scheme $X$ defined over a field $k$ with \'etale
fundamental group $\pi_1^{\et}(X)$, with $G$ the image of $\rho$, they
consider the structure (e.g., finiteness properties) of sets of $x\in
X(k_1)$, for some $k_1/k$, such that the image in $G$ of the natural
map
$$
\Gal(\bar{k}_1/k_1)\ra G
$$
associated to $x$ is not \emph{open} (in the $\ell$-adic topology), or
has large codimension in $G$, etc. In~\cite[Th. 1.1]{ct1}, general
conditions on $\rho$ are found which imply that imply that such sets
of $x\in X(k)$ are finite when $X/k$ is a smooth curve and $k$
finitely generated over $\Q$ and in~\cite[Th. 1.1]{ct2}, this is
extended to all $x\in X(k_1)$ with $[k_1:k]\leq d$ for $d\geq 1$. The
strategy parallels ours: a suitable \emph{tower} $(X_{n+1}\ra
X_{n})_{n\geq 0}$ of coverings is constructed so that its rational
points control the desired set, and Cadoret and Tamagawa show that
either the genus~\cite{ct1} or the gonality~\cite{ct2} of the curves
$X_n$ tends to infinity. However, the details of the proofs of this
facts are strikingly different from our
Theorem~\ref{th:gonalitymain}. 
\par
Note that although we have not considered such ``vertical'' towers of
coverings in this paper, our results can also be applied in this
context. Indeed, works of Bourgain and Gamburd~\cite{bg1,bg2} show
that families of Cayley graphs of $\SL_{d}(\Z/p^n\Z)$ are expanders,
for $d$ and $p$ fixed and $n$ varying, and recently Dinai~\cite{dinai,
  dinai2} has proved the polylogarithmic growth of the diameter
(implying the esperantist property by~(\ref{eq-dsc})) for families
$(\barG(\Z/p^n\Z))_{n\geq 1}$ wen $\barG/\Z_p$ is an arbitrary split
semisimple algebraic group. The argument is very elementary, and leads
to explicit constants $(c,A)$ in~(\ref{eq:esperanto}), which is of
great interest when thinking of effectiveness (as discussed below).

\subsection{Higher-dimensional families}

Our work is intrinsically limited to one-parameter families, in at
least two ways: (1) the use of gonality of curves to deduce
diophantine consequences through the theorem of Faltings; (2) the use
of the Li-Yau inequality to relate gonality to the Laplace operator
and then the combinatorial laplacian of graphs. It would be quite
interesting to know whether any similar result holds when dealing with
families of coverings of higher-dimensional varieties when the
associated Cayley-Schreier graphs are expanders or esperantist. 

\subsection{Extension to positive characteristic}

A basic question suggested by Theorem~\ref{th:main} is the following:
what happens when the base field $k$ is a global field of positive
characteristic? It is easy to modify the assumption so that it makes
sense, and we know that some version of the first part of
Theorem~\ref{th:gonalitymain} extends (as follows from~\cite[Prop. 5,
Prop. 7]{eehk}). However, the crucial step where we use the theorem of
Li and Yau is not available for the gonality argument. Moreover, a
na\"ive idea of ``lifting'' to characteristic zero (if possible) runs
into difficulties, since the gonality of a lift might be larger than
that over $k$. It would be very interesting to know if the analogue of
this gonality bound is true over all global fields. We therefore raise
the following question:

\begin{que}
  Let $k$ be a global field of positive characteristic $p>0$, or a
  finite field. Let $U/k$ be a smooth geometrically connected curve,
  and let $(U_i)$ be a family of finite, tamely ramified, \'etale
  covers of $U$, such that the Cayley-Schreier graphs of the finite
  quotient sets
$$
\pi_1^{\et}(U)/\pi_1^{\et}(U_i)
$$
(with respect to a fixed finite, symmetric, set $S$ of topological
generators of the tame fundamental group $\pi_1^{tame}(U)$) form an
expander graph.  Is it true, or not, that we necessarily have
\begin{equation}
\lim_{i\ra +\infty}{\gamma(U_i)}=+\infty?
\label{gonalitygrows}
\end{equation}
\end{que}

Note that in this setting we do have
$$
\lim_{i\ra +\infty}{g(U_i)}=+\infty
$$
since in the tamely ramified case we can lift the covering $U_i \ra U$
to a field $K$ of characteristic $0$ without changing the genus of
either curve, at which point we can embed $K$ in $\C$ and use the
arguments of the present paper.  (To be precise, we should require
that $S$ is the image in $\pi_1^{\et}(U)$ of some generating set of the
\emph{discrete} group $\pi_1(U_{\C})$.)

Poonen has shown~\cite{poon:gonality} that \eqref{gonalitygrows} holds
when $U = X(1)$ is the moduli space of elliptic curves and $(U_i)$ is
a sequence of modular curves of increasing level.

It is unclear to us whether any results along the lines of those
proved here can be expected when wild ramification is allowed.  As a
cautionary note we remark that, by contrast with the present paper,
Abyankhar has constructed many (wildly ramified) coverings $U_i \ra
\A^1/\F_{p^e}$ whose Galois groups are linear groups over fields of
characteristic $p$, and where $U_i$ has genus $0$.  On the other hand,
in the contexts considered here (covers coming from $\ell$-torsion
points of an abelian scheme $\mathcal{A}$ over $U$, with $\ell$ large
relative to the other data) these pathologies can perhaps be avoided;
in our earlier paper \cite{eehk} we show that for some families of
covers of this kind one can indeed show that $g(U_i)$ is unbounded.

\subsection{Issues of effectivity}\label{ssec-effective}

Because of its dependency on Faltings's theorem, there is currently no
chance of being able to effectively compute sets like
$$
\bigcup_{[k_1:k]=d}{U_i(k_1)},
$$
with notation as in Theorem~\ref{th:main}, even if we know that it is
a finite set. 
\par
However, a more accessible kind of effectivity would be to ask for an
effective determination, for a fixed $d\geq 1$, of a finite set
$I_d\subset I$ of exceptional $i$ such that the set above is finite
when $i\notin I_d$. In the context of Theorem~\ref{th:torsion}, for
instance, this would mean finding an effective $\ell_0=\ell_0(d)$ such
that
$$
\bigcup_{[k_1:k]=d}{\{t\in U(k_1)\,\mid\,
  \text{$\mathcal{A}_t[\ell](k_1)$ is non-zero} \}} 
$$
is finite if $\ell\geq \ell_0$. Or one might ask for an effective
growth result for gonality in the context of
Theorem~\ref{th:gonalitymain}.
\par
As our argument shows, this is directly related to the issue of
finding effective expansion constants for the families of Cayley
graphs that we use,\footnote{\ The other ingredient which we can not
  obviously estimate is the first Neumann eigenvalue $\eta$ that
  occurs in the proof of the comparison inequality~(\ref{eq-burger}),
  as described in Appendix A below (see~\ref{eq:mu}).} which is a
delicate question in general, even if the methods of~\cite{bgt}
and~\cite{psz} are effective in principle. Here using vertical towers
can lead to drastic simplifications, as shown by the results of
Dinai~\cite{dinai, dinai2}, where simple explicit constants $(c,A)$
are obtained. We hope to use these to make progress towards
effectivity of some our results.

\section*{Appendix A: the Burger method}
\label{app:burger}

In this appendix, we sketch the extension of Burger's comparison
principle to finite-area hyperbolic surfaces, as required in our
arguments. The arguments follow Burger's method (most clearly
explained in his thesis~\cite[Ch. 6]{burger3}, which is not readily
available, and only briefly sketched in~\cite{burger,burger2}).

\begin{thm}[Burger]\label{th:burger}
  Let $U'\ra U$ be a finite covering of a connected hyperbolic Riemann
  surface with finite hyperbolic area. Fix a symmetric system of
  generators $S$ of $\pi_1(U,x_0)$ and let
$$
\Gamma=C(\pi_1(U,x_0)/\pi_1(U',x'_0),S)
$$
be the associated Cayley-Schreier graph, where $x'_0\in U'$ is a point
above $x_0$. Then there exists a constant $c>0$, depending only on $U$
and $S$, such that
$$
\lambda_1(U')\geq c\lambda_1(\Gamma).
$$
\end{thm}

Let $P$ be the set of points in $U'$ above $x_0$, $\tilde{P}$ the set
of those in the universal cover $\tilde{U}=\mathbb{H}$ of $U$, so that
$|P|=|\Gamma|$. For $x\in P$, $\tilde{x}\in\tilde{P}$, let
\begin{gather*}
  \mathcal{F}(x)=\{u\in U'\,\mid\, d(u,x)<d(u,x')\text{ for all }
  x'\in
  P,\ x'\not=x\},\\
  \tilde{\mathcal{F}}(\tilde{x})=\{u\in \tilde{U}\,\mid\,
  d(u,\tilde{x})<d(u,x')\text{ for all } x'\in \tilde{P},\
  x'\not=\tilde{x}\}.
\end{gather*}
\par
It is well-known that each $\tilde{\mathcal{F}}(\tilde{x})\subset
\tilde{U}$ is a fundamental domain for the action of $\pi_1(U,x_0)$ on
$\tilde{U}$, and $\mathcal{F}(x)\subset U'$ is one for the covering
$U'\ra U$. When $\tilde{x}$ (resp. $x$) varies, these are
disjoint. The closures $\overline{\mathcal{F}(x)}$ cover $U'$, with
boundaries having measure $0$, and hence
$$
\mu(\mathcal{F}(x))=\mu(U)<+\infty
$$
for all $x\in P$. Moreover, the set
$$
T=\{g\in\pi_1(U,x_0)\,\mid\, g\not=1,\
g\overline{\tilde{\mathcal{F}}(\tilde{x})} \cap
\overline{\tilde{\mathcal{F}}(\tilde{x})}\not=\emptyset\}
$$ 
is a finite symmetric generating set of $\pi_1(U,x_0)$. It is an
elementary fact that we need only prove the result when the generating
set $S$ is replaced by $T$. We denote
$$
r=|T|+1.
$$
\par
Now consider the graph $\Gamma'$ with vertex set $P$ and edges joining
$x$ and $x'$ in $P$, $x\not=x'$, if and only if
$$
\overline{\mathcal{F}(x)}\cap \overline{\mathcal{F}(x')}\not=\emptyset.
$$
\par
This graph may be non-regular, but its valence function $v\,:\, P\ra
\R$ satisfies $1\leq v(x)\leq |T|$ for all $x$. In fact, one checks
that $\Gamma'$ is obtained from the Cayley-Schreier graph
$$
\Gamma_T=C(\pi_1(U,x_0)/\pi_1(U',x'_0),T)
$$ 
by (i) replacing multiple edges by simple ones; (ii) removing
loops. For simplicity, we will assume that $\Gamma'=\Gamma_T$, and in
particular that $v(x)=|T|$ for all $x$. (See also~\cite[Ch.3, \S
4]{burger3} for details about this construction.)
\par
To prove the theorem, we use the variational characterization (or
definition, see~(\ref{eq:lambdaone})) of $\lambda_1=\lambda_1(U')$
$$
\lambda_1=\inf\Bigl\{ \frac{\displaystyle{\int_{U_i}{\|\nabla
        \varphi\|^2d\mu}}}{\|\varphi\|^2}
\,\mid\, \varphi\text{ smooth and } \int_{U'}{\varphi(x)d\mu(x)}=0
\Bigr\},
$$
where again $\nabla\varphi$ refers to the hyperbolic gradient of
$\varphi$.
\par
Precisely, we have already recalled that either this quantity is
$=1/4$ (in which case we are done) or else there exists a non-zero
(eigenfunction) $\psi\in L^2(U',d\mu)$ with mean zero over $U'$ and
which attains the infimum. Of course, we now consider this case, and
we may assume that $\psi$ has $L^2$-norm equal to $1$.
% $$
% \Delta \psi=\lambda_1 \psi.
% $$
\par
Let $L^2(P)$ be the space of functions on $P$, with the inner product
$$
\langle g_1,g_2\rangle = \sum_{x\in P}{g_1(x)\overline{g_2(x)}}.
$$
\par
Now we must perform a transfer of some kind from smooth functions on
$U'$ to discrete functions on the vertex set $P$. The idea is quite
simple: to a function $f$, we associate the function
$$
\Phi(f)\,:\,x\mapsto \int_{\mathcal{F}(x)}{fd\mu}
$$
on $P$. This linear map $\Phi$ is of course not an isometry from
$L^2(U',d\mu)$ to $L^2(P)$, but it is at least continuous since
\begin{equation}\label{eq:cont}
\|\Phi(f)\|^2=\sum_{x\in P}{\Bigl|
\int_{\mathcal{F}(x)}{fd\mu}\Bigr|^2}
\leq \mu(U)\|f\|^2
\end{equation}
by the Cauchy-Schwarz inequality and the fact that the
$\mathcal{F}(x)$ are disjoint and have measure $\mu(U)$.
\par

\begin{proof}[Proof of Theorem~\ref{th:burger}]
  For $x\in P$, let $N(x)=\{x\}\cup \{x'\text{ adjacent to } x\}$, and
  define
$$
\mathcal{G}(x)=\bigcup_{x'\in N(x)}{\overline{\mathcal{F}(x)}}\subset
U',
$$
so that (under our assumption $\Gamma_T=\Gamma'$, and recalling that
$r=|T|+1$) we have $|N(x)|=r$ and
$$
\mu(\mathcal{G}(x))=r\mu(U).
$$
\par
We start by stating the following fact, to be proved below (this is
where the distinction between compact and finite-area surfaces will
occur):
\par
\medskip \textbf{Fact 1.} There exists a constant $\eta>0$, depending
only on $U$, such that, for all $x\in P$ and for
$\mathcal{H}=\mathcal{F}(x)$ or $\mathcal{G}(x)$, we have
\begin{equation}\label{eq:mu}
  \inf\Bigl\{
  \frac{\displaystyle{\int_{\mathcal{H}}{\|\nabla \varphi\|^2d\mu}}}{
\int_{\mathcal{H}} |\varphi|^2 d\mu}
  \,\mid\, 0\not=\varphi\text{ smooth and } 
\int_{\mathcal{H}}{\varphi(x)d\mu(x)}=0
  \Bigr\}\geq \eta.
\end{equation}
\par
\medskip
Assuming this, consider a non-zero function of the type
$$
f=\alpha+\beta \psi
$$
on $U'$ with $\alpha$, $\beta\in\R$, and $\psi$ the eigenfunction
described above.  We have an obvious inequality
\begin{equation}\label{eq:tocompare}
\int_{U'}{\|\nabla f\|^2d\mu}=\lambda_1\beta^2\leq \lambda_1\|f\|^2.
\end{equation}
\par
Now we proceed to bound the left-hand side from below using the
pieces $\mathcal{G}(x)$. For any $x\in P$, the function
$$
\varphi=f-\frac{1}{\mu(\mathcal{G}(x))}\int_{\mathcal{G}(x)}{fd\mu}
$$
(with $\nabla \varphi=\nabla f$) can be used to test~(\ref{eq:mu}),
and therefore we have
\begin{align*}
\int_{\mathcal{G}(x)}{\|\nabla f\|^2}&\geq \eta
\int_{\mathcal{G}(x)}{
\Bigl(f-\frac{1}{\mu(\mathcal{G}(x))}\int_{\mathcal{G}(x)}{fd\mu}
\Bigr)^2d\mu
}\\
&=\eta
\Bigl\{\int_{\mathcal{G}(x)}{f^2d\mu}-\frac{1}{\mu(\mathcal{G}(x))}
\Bigl(\int_{\mathcal{G}(x)}{fd\mu}\Bigr)^2\Bigr\}.
\end{align*}
\par
Since we assumed that $\Gamma'$ is regular, each $\mathcal{G}(x)$ is
the union of $r$ among the $\mathcal{F}(x')$. Therefore, if we sum
over $x\in P$, divide by $r$ and use the fact that
$\mu(\mathcal{G}(x))=r\mu(U)$, we obtain
$$
\|\nabla f\|^2\geq \eta \Bigl\{
\|f\|^2-
\frac{1}{r^2\mu(U)}\sum_{x\in P}
\Bigl(\sum_{x'\in N(x)}{\Phi(f)(x')}\Bigr)^2
\Bigr\}.
$$
\par
Comparing with~(\ref{eq:tocompare}) and dividing by $\|f\|^2$, we find
using~(\ref{eq:cont}) that
$$
\lambda_1\geq \eta \frac{\langle
  B\Phi(f),\Phi(f)\rangle}{\|\Phi(f)\|^2}
$$
where the linear operator $B$ on $L^2(P)$ is defined by
$$
B=1-\frac{1}{r^2}A^2,
$$
with $A$ being the self-adjoint linear map on $L^2(P)$ defined by
$$
A(g)(x)=\sum_{x'\in N(x)}{g(x)}.
$$
\par
The crucial point is that since the combinatorial Laplace operator
$\Delta$ of $\Gamma'$ is given by $\Delta=r\mathrm{Id}-A$ (where the
assumption $\Gamma'=\Gamma_T$ is used again), the operator $B$ is
itself closely related to $\Delta$, namely
\begin{equation}\label{eq:b}
B=\frac{1}{r^2}(r^2\mathrm{Id}-A^2)=\frac{1}{r^2}\Delta (2r-\Delta).
\end{equation}
\par
It is clear that $B$ is $\geq 0$ and has eigenvalue $0$ with
multiplicity $1$ for the constant eigenfunction. Let $\lambda'_1>0$
denote the smallest positive eigenvalue of $B$. We claim:
\par
\medskip \textbf{Fact 2.} There exists $c>0$, depending only on $U$,
such that either $\lambda_1\geq c$ or else there exists $\alpha$,
$\beta\in\R$ not both zero for which $\Phi(f)=\Phi(\alpha+\beta\psi)$
is non-zero and has mean zero on $P$.
\par
\medskip
If this is the case, we construct the test function $f$ using these
$\alpha$ and $\beta$; then, since 
$$
\langle \Phi(f),1\rangle =0,
$$
by the variational inequality for the spectrum of $B$, we have
$$
\frac{\langle B\Phi(f),\Phi(f)\rangle}{\|\Phi(f)\|^2}\geq \lambda'_1. 
$$
\par
However, using~(\ref{eq:b}), we can compare $\lambda'_1$ and the first
eigenvalue $\lambda_1(\Gamma')$: from $\|\Delta g\|\leq r\|g\|$ for
all $g$, we get
$$
\lambda'_1\geq
\frac{1}{r}\lambda_1(\Gamma')=\frac{1}{r}\lambda_1(\Gamma_T)
$$
by looking on the subspace 
$$
L^2_0(P)=(\C\cdot 1)^{\perp}\subset L^2(P),
$$ 
which is stable under $B$, $\Delta$ and $2r-\Delta$, and where each of
these operators is invertible: indeed, we have
$$
\frac{1}{\lambda'_1}=\|B^{-1}\|\leq
r^2\|\Delta^{-1}\|\|(2r-\Delta)^{-1}\| \leq
r\|\Delta^{-1}\|=\frac{r}{\lambda_1(\Gamma')},
$$
all operators and norms thereof being computed on $L^2_0(P)$
(see~\cite[p. 73, (d), (e)]{burger3} or~\cite[\S 3, Cor. 1
(a)]{burger2} for the general case where $\Gamma'\not=\Gamma_T$;
Burger shows that $\lambda'_1\geq \frac{2}{r^3}\lambda_1(\Gamma)$).
\par
Combining these inequalities, we find that
$$
\lambda_1\geq \min\Bigl(c,\frac{\eta}{r}\lambda_1(\Gamma_T)\Bigr),
$$
which concludes our proof.
\par
We now justify the two claims above.  For Fact 1, we note that the
infimum considered, say $\eta(\mathcal{H})$, are nothing but the
smallest positive eigenvalue for the Laplace operator with Neumann
boundary condition on $\mathcal{H}$ -- or more precisely, because the
area of $\mathcal{H}$ is finite, the constant function $1$ gives the
base eigenvalue $0$ as before, and $\eta(\mathcal{H})$ is either $1/4$
or the first positive eigenvalue (by~\cite[Th. 2.4]{ph-sa}, which
states in much greater generality that the spectrum is discrete in
$[0,1/4]$; if $\mathcal{H}$ were compact, this becomes standard
spectral geometry of compact Riemannian manifolds). So
$\eta(\mathcal{H})>0$, but we must still show that there is a
lower-bound depending only on $U$, not on $U'$.
\par
For this, fix any $\tilde{x}\in \tilde{U}$ above $x_0$. The reasoning
in~\cite[p. 71]{burger3} applies identically to show that
$\mathcal{G}(x)$ is always isometric to a quotient of the domain
$$
\mathcal{A}= \bigcup_{s\in T\cup
  \{1\}}{s\tilde{\mathcal{F}}(\tilde{x}_1)}\subset \tilde{U}=\mathbb{H},
$$
(which depends only on $U$) under an equivalence relation of
congruence modulo a subset of $T^3$. Each such quotient is a
finite-area domain in $\mathbb{H}$, hence also its first non-zero
Neumann eigenvalue (defined variationally) is $>0$
by~\cite[Th. 2.4]{ph-sa}. Since $T$ is finite, there are only finitely
many such quotients to consider, depending only on $U$, hence the
smallest among these Neumann eigenvalues is still $\eta>0$, and we
have of course
$$
\eta(\mathcal{H})\geq \eta
$$
for all $\mathcal{H}$, proving Fact 1.
\par
For Fact 2, we note that to find the required test function $f$ it is
enough to know that the map $\Phi$ is injective on the $\R$-span of
$1$ and $\psi$. Indeed, we can then find a non-zero $f$ in the kernel
of the linear functional
$$
f\mapsto \langle \Phi(f),1\rangle
$$
(which is non-trivial since $1$ maps to $|P|$), and this $f$ will
satisfy the required conditions.
\par
Now we have two cases. If $\lambda_1\geq \eta$, where $\eta$ is given
by Fact (1), we are done (and take $c=\eta$). Otherwise, we have
\begin{equation}\label{eq:assump}
0<\lambda_1<\eta,
\end{equation}
and we now show that this implies that $\Phi$ is injective on the
(real) span of $1$ and $\psi$, which thus concludes the proof.
\par
Thus, let $\alpha$, $\beta\in\R$ be such that
$\Phi(f)=\Phi(\alpha+\beta \psi)=0$. Then, for all $x\in P$, we have
$$
\int_{\mathcal{F}(x)}{\|\nabla f\|^2d\mu}\geq \eta
\int_{\mathcal{F}(x)}{f^2d\mu},
$$
since $\Phi(f)=0$ means that $f$ restricted to $\mathcal{F}(x)$ can be
used to test~(\ref{eq:mu}). Summing over $x$, we get
$$
\|\nabla f\|^2\geq \eta \|f\|^2,
$$
but $f=\alpha+\beta \psi$ implies then that
$$
\eta\|f\|^2\leq \|\nabla f\|^2=\beta^2\lambda_1\leq \lambda_1\|f\|^2,
$$
and by comparing with~(\ref{eq:assump}), we see that $f=0$.
\end{proof}

\section*{Appendix B: Semisimple Approximation \`a la Nori and Serre}
\label{app:noriserre}

\setcounter{ccounter}{0}

For $m$ a positive integer and for $\ell$ varying over the primes,
there are two kinds of groups that we focus on in this section.  The
first are the (finite) semisimple subgroups $G\leq\GL_m(\F_\ell)$,
that is, subgroups which act semisimply in the natural representation
of $\GL_m(\F_\ell)$ on $V=\F_\ell^m$, and the second are connected
semisimple groups $\barG\subseteq\barGL_m/\F_\ell$.  Nori showed that,
outside of an explicit finite set of exceptional $\ell$, there is a
bijection between the finite semisimple $G$ which are generated by
their elements of order $\ell$ and ``exponentially-generated'' $\barG$
(see \cite{nori:inv}).  Serre showed that, up to excluding finitely
many more $\ell$, the semisimple groups $\barG$ which occur come from
a finite collection of groups in characteristic zero
(cf.~\cite{serre:vig}).  We will give a brief review of these results
following both \cite{serre:vig} and a set of notes taken during the
course mentioned in \opcit.
\par
Let $\ell$ be a prime and $G\leq\GL_m(\F_\ell)$ be a (finite)
subgroup.  We write $G_u\subseteq G$ for the subset of unipotent
elements and $G^+\leq G$ for the characteristic subgroup
$G^+=\gp{G_u}$.

\begin{prop}\label{prop:order}
If $\ell\geq m-1$ and if $G\leq\GL_m(\F_\ell)$, then $g^\ell=1$,
for every $g\in G_u$, and $G/G^+$ has order prime to $\ell$.
\end{prop}

\begin{proof}
  Let $P\leq G$ be an $\ell$-Sylow subgroup of $G$.  Every element
  $g\in P$ lies in $G_u$, and thus $(g-1)^{m-1}=0$.  On the other
  hand, because $\ell\geq m-1$, $g^\ell-1=(g-1)^\ell=0$, thus $g$ is
  killed by $\ell$.  Therefore $P\leq G^+$ and $G/G^+$ has
  order prime to $\ell$ because $G/P$ does.
\end{proof}

If $\ell\geq m$, the exponential and logarithm maps give
mutually-inverse bijections between the unipotent elements of
$\GL_m(\F_\ell)$ and the nilpotent elements of $\M_m(\F_\ell)$, and we
write $\g\leq\M_m(\F_\ell)$ for the $\F_\ell$-span of $\log(G_u)$.

\begin{prop}\label{prop:restricted}
  If $\ell\geq 2m-1$, then $\g$ is an $\F_\ell$-Lie subalgebra of
  $\M_m(\F_\ell)$ and every $\g$-submodule of
  $\bar{V}=V\otimes\Fbar_\ell$ is a $G^+$-submodule.
\end{prop}

\begin{proof}
  This follows from Lemmas 1.4 (applied with $W_1=W_2=W$) and 1.6 of
  \cite{nori:inv}.
\end{proof}

For each $g\in G_u-\{1\}$, we can use the embedding
$\gp{\log(g)}\to\M_m(\F_\ell)$ to extend the embedding
$\gp{g}\to\GL_m(\F_\ell)$, to a one-parameter subgroup
$\G_a\to\barGL_m$ over $\F_\ell$.  We write
$\barG^+\subset\barGL_m$ for the algebraic subgroup generated by the
images of all such one-parameter subgroups and
$\Lie(\barG^+)\subseteq\M_m(\F_\ell)$ for the $\F_\ell$-Lie subalgebra
of $\barG^+\subset\barGL_m$.

\begin{thm}\label{thm:noriB}
  There is a constant $\ell_1=\ell_1(m)\geq 2m-1$ such that if
  $\ell\geq\ell_1$ and if $G\leq\GL_m(\F_\ell)$, then
  $G^+=\barG^+(\F_\ell)^+$ and $\g=\Lie(\barG^+)$.
\end{thm}

%Note, for each $g\in G_u-\{1\}$, the corresponding one-parameter subgroup
%$\underline{\G}_a\subset\barGL_m$ lies in $\barG^+\subseteq\barGL_m$ by
%definition, hence $\log(g)\in\Lie(\underline{\G}_a)\subseteq\Lie(\barG^+)$.
%Thus there is an a priori inclusion $\g\subseteq\Lie(\barG^+)$.

\begin{proof}
  This is Theorem B of \cite{nori:inv}.
\end{proof}

We call the rational representation $\barG^+\to\barGL_m$ the algebraic
envelope of $G^+\to\GL_m(\F_\ell)$.  We are most interested in the
case where $G$ acts semisimply, and then the following corollary shows
that, for almost all $\ell$, the algebraic envelope
$\barG^+\to\barGL_m$ is a rational representation of a semisimple
group.  We will see that there are strong restrictions on the dominant
weights occurring in this representation and that there are finitely
many $\Z$-groups which give rise to them.

\begin{cor}\label{cor:noriB}
  If $\ell\geq\ell_1$ and if $G\leq\GL_m(\F_\ell)$ is semisimple, then
  $\barG^+$ is semisimple.
\end{cor}

\begin{proof}
  By assumption, $G$ acts semisimply on $V$, so Clifford's Theorem
  (see, e.g.,~\cite[Theorem 49.2]{cr}) implies that
  $G^+=\barG^+(\F_\ell)^+$ acts semisimply on $V$.  Because $\barG^+$
  is exponentially generated, the radical of $\barG^+$ is unipotent,
  so we denote it $\barU$.  % Clifford's Theorem (see,
  % e.g.,~\cite[Theorem 49.2]{cr}) implies $\barG^+(\F_\ell)^+=G^+$ acts
  % semisimply on $V$, and 
  Another application of Clifford's Theorem implies that
  $\barU(\F_\ell)$ also acts semisimply on $V$, and hence is trivial.
  But $\barU(\F_\ell)$ has $\ell^{\dim(\barU)}$ elements, and so the
  triviality of the group $\bar{U}(\F_\ell)$ implies that $\barU$
  itself is trivial as algebraic group. Thus $\barG^+$ is semisimple.
\end{proof}

For the remainder of this section we suppose $\ell\geq\ell_1$ and
$G\leq\GL_m(\F_\ell)$ acts semisimply and satisfies $G=G^+$, i.e., it
is generated by elements of order $\ell$.  We write $\barG\to\barG^+$
for the simply-connected cover of $\barG$, and $\barG\to\barGL_m$ for
the induced rational representation.  The hypotheses on $\ell$ and $G$
imply that $\barG/\Fbar_\ell$ is a simply-connected semisimple group
of rank at most $m-1$, the rank of $\barSL_m$.  Moreover, there is a
finite collection $\{\barG_i\to\Spec(\Z)\}$ of split simply-connected
semisimple groups (independent of $\ell$ and $G$) such that, for some
$i$, the group $\barG/\Fbar_\ell$ is isomorphic to $\barG_i/\Fbar_\ell$.

For each $i$, the group $\barG_i\to\Spec(\Z)$ is a simply-connected
Chevalley group.  If we fix a maximal torus $\barT_i\subset\barG_i$
over $\Z$, then the irreducible representations of $\barG_i/\C$ are
parametrized by their dominant weights $\lambda\in X(\barT_i)_+$.
Steinberg showed that there are $\Z$-forms
$\rho_\lambda:\barG_i\to\barGL(V_\lambda)$ of these representations
(see \cite{steinberg}), and one can show there is an explicit constant
$\ell(\lambda)$ such that, for every $\ell\geq\ell(\lambda)$, the
fiber $\rho_\lambda/\Fbar_\ell$ is also irreducible.  If we fix a set
$\{\w_{ij}\}$ of fundamental weights of $\barT_i$, then one can also
show that the finite subset $\Lambda_i\subset X(\barT_i)_+$ of
dominant weights $\lambda=\Sigma_j c_j\w_{ij}$ satisfying
$\max\{c_j\}\leq m-1$ contains all $\lambda$ satisfying
$\dim(V_\lambda)\leq m$.  We will see that, for almost all $\ell$,
each irreducible subrepresentation of $\barG\to\barGL_m$ over
$\Fbar_\ell$ is isomorphic to $\rho_\lambda/\Fbar_\ell$ for some
$\rho_\lambda:\barG_i\to\barGL(V_\lambda)$ with $\lambda\in\Lambda_i$.

We fix a maximal torus $\barT\subset\barG$ and a set $\{\w_i\}$ of
fundamental weights of $\barT$, and let $\Lambda_\ell\subset X(\barT)$
be the finite set of dominant weights $\lambda=\Sigma_i c_i\w_i$ which
occur in $\barG\to\barGL_m$.  The following proposition shows that the
weights $\lambda\in\Lambda_\ell$ are $\ell$-restricted, and thus the
rational representation $\barG\to\barGL_m$ is restricted.

\begin{prop}
  If $\lambda=\Sigma_i c_i\w_i\in\Lambda_\ell$, then $c_i\leq\ell-1$.
\end{prop}

\begin{proof}
  On the one hand, the irreducible submodules for $G$ and $\barG$ of
  $\bar{V}=V\otimes\Fbar_\ell$ coincide, and
  Proposition~\ref{prop:restricted} implies they are all
  $\g$-irreducible.  On the other hand, the only irreducible
  $\barG$-modules over $\Fbar_\ell$ which are $\g$-irreducible are
  those whose dominant weight $\lambda=\Sigma_i c_i\w_i$ satisfies
  $c_i\leq\ell-1$ (cf.~\cite[Part II, Section 3.15]{jantzen}).
\end{proof}

A priori the set $\Lambda_\ell$ could grow with $\ell$, but the
following proposition shows that it is bounded in a very strong sense.

\begin{prop}
  If $\lambda=\Sigma_i c_i\w_i\in\Lambda_\ell$, then $c_i\leq m-1$.
\end{prop}

\begin{proof}
  By the previous proposition, $\lambda$ is $\ell$-restricted.  On the
  one hand, for $n=c_i$, the rational representation
  $\barSL_2\to\barGL_{n+1}$ corresponding to the $\ell$-restricted
  weight $c_i\w_i$ is the $n$-th symmetric power of the standard
  representation $\barSL_2\to\barGL_2$ and is irreducible
  (cf.~\cite[Part II, Section 3.0]{jantzen}).  On the other hand, for
  each $i$, there is an embedding $\barSL_2\to\barG$ such that
  $\lambda_i=c_i\w_i$ is one of the dominant weights of the induced
  representation $\barSL_2\to\barGL_m$, thus $m\geq n+1=c_i+1$.
\end{proof}

We already saw that $\barG/\Fbar_\ell$ is isomorphic to
$\barG_i/\Fbar_\ell$, for some $i$, and together with the last
proposition we complete the proof of the claim that the dominant
weights $\lambda$ occurring in $\barG\to\barGL_m$ lie in $\Lambda_i$.
The upshot is that we obtain the following theorem.

\begin{thm}
  There exists a finite collection $\{\rho_{ij}:\barG_i\to\barGL_m\}$
  of $\Z$-representations of simply-connected Chevalley groups and a
  constant $\ell_2=\ell_2(m)\geq\ell_1$ such that if $\ell\geq\ell_2$
  and if $G\leq\GL_m(\F_\ell)$ is semisimple and satisfies $G=G^+$,
  then for some $i,j$, the fiber $\rho_{ij}/\Fbar_\ell$ is isomorphic
  to $\barG\to\barGL_m$.
\end{thm}

For each pair of integers $r,s\geq 1$, we write $T_{rs}V$ for the
vector space $T_{rs}V = (\oplus_{i=1}^{r} V^{\otimes i})^{\oplus s}$
and $\barGL(V)\to\barGL(T_{rs}V)$ for the corresponding tensor
representation.

\begin{cor}\label{cor:tensor}
  There are constants $\ell_3=\ell_3(m)\geq\ell_2$, $r=r_1(m)$, and
  $s=s_1(m)$ such that if $\ell\geq\ell_3$ and if
  $G\leq\GL_m(\F_\ell)$ is semisimple and satisfies $G=G^+$, then the
  composite representation $\barG\to\barGL(V)\to\barGL(T_{rs}V)$
  identifies $\barG$ with the algebraic subgroup of elements in
  $\barGL(V)$ acting trivially on the subspace of $\barG$-invariants
  in $T_{rs}V$.
\end{cor}

\begin{proof}
  The main idea is to show that, for each $i,j$, an analogous
  statement holds for the $\Q$-fiber of
  $\rho_{ij}:\barG_i\to\barGL_m$.  More precisely, if we write
  $V_\Q=\Q^m$, then for some $r=r(i,j,m)$ and $s=s(i,j,m)$ depending
  on $i$ and $j$, the tensor representation $\barG_i(\Q)\to\GL(T_{rs}
  V_\Q)$ identifies $\barG_i/\Q$ with the algebraic subgroup of
  $\barGL_m$ acting trivially on the subspace of
  $\barG_i$-invariants.\footnote{.\ Every irreducible
    finite-dimensional rational representation of $\barGL(V_\Q)$ is a
    subquotient of $\barGL(V_\Q)\to\barGL(V_\Q^{\otimes r})$ for some
    $r$, so every finite-dimensional rational representation is a
    subquotient of $\barGL(V_\Q)\to\barGL(T_{rs} V_\Q)$ for some
    $r,s$.  In particular, there is a finite-dimensional rational
    representation $\rho:\barGL(V_\Q)\to\barGL(W_\Q)$ which identifies
    $\barG_i$ with the stabilizer in $\barGL(V_\Q)$ of a line in
    $W_\Q$ (see, e.g.,~Corollary 3.5 of \cite[II.2]{dg}), and because
    $\barG_i$ is semisimple, it acts trivially on the line.}
  Moreover, there is a $\Z$-form of this tensor representation, and
  for almost all $\ell$, the corresponding tensor representation
  $\barG_i(\F_\ell)\to\GL(T_{rs} V)$ identifies $\barG_i/\F_\ell$ with
  the algebraic subgroup of $\barGL(T_{rs} V)$ acting trivially on the
  subspace of $\barG_i$-invariants.  In particular, in light of the
  theorem it suffices to take $r_1(m)\geq\max\{r(i,j,m)\}$ and
  $s_1(m)\geq\max\{s(i,j,m)\}$.
\end{proof}

%%%%%%%%%%%%%%%%%%%%%%%%%%%%%%%%%%%%%%%%%%%%%%%%%%%%%%%%%%%%%%%%%%%%%%%%

\end{document}